\documentclass{amsart}
\usepackage{amssymb,amsmath,stmaryrd,mathrsfs,amsthm}
\usepackage[all,2cell]{xy}
\usepackage[neveradjust]{paralist}
\usepackage{hyperref}

\def\mdef#1#2{\expandafter\expandafter\expandafter\gdef\expandafter\expandafter\noexpand#1\expandafter{\expandafter\ensuremath\expandafter{#2}}}

\newcommand{\cP}{\ensuremath{\mathcal{P}}}

\newcommand{\cS}{\ensuremath{\mathcal{S}}}
\newcommand{\cT}{\ensuremath{\mathcal{T}}}

\newcommand{\bbN}{\ensuremath{\mathbb{N}}}

\newcommand{\bbQ}{\ensuremath{\mathbb{Q}}}
\newcommand{\bbR}{\ensuremath{\mathbb{R}}}
\newcommand{\bbS}{\ensuremath{\mathbb{S}}}

\newcommand{\bbZ}{\ensuremath{\mathbb{Z}}}

\DeclareSymbolFont{bbold}{U}{bbold}{m}{n}
\DeclareSymbolFontAlphabet{\mathbbb}{bbold}

\newcommand{\bbk}{\ensuremath{\mathbbb{k}}}

\newcommand{\bA}{\ensuremath{\mathbf{A}}}
\newcommand{\bB}{\ensuremath{\mathbf{B}}}
\newcommand{\bC}{\ensuremath{\mathbf{C}}}

\newcommand{\bS}{\ensuremath{\mathbf{S}}}

\newcommand{\bV}{\ensuremath{\mathbf{V}}}
\newcommand{\bW}{\ensuremath{\mathbf{W}}}

\newcommand{\ten}{\ensuremath{\otimes}}

\newcommand{\inv}{^{-1}}

\newcommand{\Id}{\ensuremath{\operatorname{Id}}}

\newcommand{\colim}{\ensuremath{\operatorname{colim}}}

\newcommand{\op}{\ensuremath{^{\mathit{op}}}}

\newdir{ >}{{}*!/-10pt/@{>}}
\newcommand{\iso}{\cong}

\newcommand{\too}[1][]{\ensuremath{\overset{#1}{\longrightarrow}}}
\newcommand{\oot}[1][]{\ensuremath{\overset{#1}{\longleftarrow}}}

\newcommand{\toto}{\ensuremath{\rightrightarrows}}
\newcommand{\into}{\ensuremath{\hookrightarrow}}

\newcommand{\mono}{\rightarrowtail}

\newcommand{\maps}{\colon}

\let\xto\xrightarrow

\newlength{\htolgth}

\newtheorem{thm}{Theorem}[section]
  
\newtheorem{cor}{Corollary}
  \makeatletter\let\c@cor\c@thm\makeatother
  \numberwithin{cor}{section}

  \makeatletter\let\c@prop\c@thm\makeatother
  \numberwithin{prop}{section}

  \makeatletter\let\c@lem\c@thm\makeatother
  \numberwithin{lem}{section}

  \makeatletter\let\c@sch\c@thm\makeatother
  \numberwithin{sch}{section}

  \makeatletter\let\c@assume\c@thm\makeatother
  \numberwithin{assume}{section}

  \makeatletter\let\c@claim\c@thm\makeatother
  \numberwithin{claim}{section}

  \makeatletter\let\c@conj\c@thm\makeatother
  \numberwithin{conj}{section}

  \makeatletter\let\c@hyp\c@thm\makeatother
  \numberwithin{hyp}{section}

  \makeatletter\let\c@fact\c@thm\makeatother
  \numberwithin{fact}{section}
  
\theoremstyle{definition}

  \makeatletter\let\c@defn\c@thm\makeatother
  \numberwithin{defn}{section}

  \makeatletter\let\c@notn\c@thm\makeatother
  \numberwithin{notn}{section}
  
\theoremstyle{remark}
\newtheorem{rmk}{Remark}
  \makeatletter\let\c@rmk\c@thm\makeatother
  \numberwithin{rmk}{section}
  
\newtheorem{eg}{Example}
  \makeatletter\let\c@eg\c@thm\makeatother
  \numberwithin{eg}{section}

  \makeatletter\let\c@egs\c@thm\makeatother
  \numberwithin{egs}{section}

  \makeatletter\let\c@ex\c@thm\makeatother
  \numberwithin{ex}{section}

  \makeatletter\let\c@ceg\c@thm\makeatother
  \numberwithin{ceg}{section}

\makeatletter
\def\alwaysmath#1{\expandafter\expandafter\expandafter\global\expandafter\expandafter\expandafter\let\expandafter\expandafter\csname your@#1\endcsname\csname #1\endcsname
  \expandafter\def\csname #1\endcsname{\ensuremath{\csname your@#1\endcsname}}}
\alwaysmath{alpha}
\alwaysmath{beta}
\alwaysmath{gamma}
\alwaysmath{Gamma}
\alwaysmath{delta}
\alwaysmath{Delta}
\alwaysmath{epsilon}

\alwaysmath{zeta}
\alwaysmath{eta}
\alwaysmath{theta}
\alwaysmath{Theta}
\alwaysmath{iota}
\alwaysmath{kappa}
\alwaysmath{lambda}
\alwaysmath{Lambda}
\alwaysmath{mu}
\alwaysmath{nu}
\alwaysmath{xi}
\alwaysmath{pi}
\alwaysmath{rho}
\alwaysmath{sigma}
\alwaysmath{Sigma}
\alwaysmath{tau}
\alwaysmath{upsilon}
\alwaysmath{Upsilon}
\alwaysmath{phi}
\alwaysmath{Phi}
\newcommand{\ph}{\ensuremath{\varphi}}
\alwaysmath{chi}
\alwaysmath{psi}
\alwaysmath{Psi}
\alwaysmath{omega}
\alwaysmath{Omega}
\alwaysmath{ell}
\alwaysmath{infty}
\makeatother

\let\al\alpha
\let\be\beta
\let\gm\gamma

\let\de\delta

\let\si\sigma

\let\om\omega
\let\ka\kappa
\let\la\lambda

\title{Set theory for category theory}
\author{Michael A. Shulman}
\date{\today}
\newcommand{\Set}{\ensuremath{\mathbf{Set}}}
\mdef\bAbar{\overline{\mathbf{A}}}
\mdef\Mnd{\mathbf{Mnd}}
\mdef\lbb\llbracket
\mdef\rbb\rrbracket
\mdef\Def{\mathrm{Def}}
\mdef\cf{\mathrm{cf}}
\let\ee\prec
\def\zfc{\textsc{zfc}}
\def\zfci{\textsc{zfc+i}}
\def\zfcs{\textsc{zfc/s}}
\def\zmc{\textsc{zmc}}
\def\zmcs{\textsc{zmc/s}}
\def\nbg{\textsc{nbg}}
\def\mk{\textsc{mk}}

\def\wptnc{\textsc{wptnc}}
\def\etcs{\textsc{etcs}}
\def\bzc{\textsc{bzc}}
\def\bzci{\textsc{bzc+i}}
\def\mac{\textsc{mac}}
\def\ch{\textsc{ch}}
\def\nfu{\textsc{nfu}}
\def\nno{\textsc{nno}}
\def\ast{\textsc{ast}}

\begin{document}
\maketitle

\begin{abstract}
  Questions of set-theoretic size play an essential role in category
  theory, especially the distinction between sets and proper classes
  (or small sets and large sets).  There are many different ways to
  formalize this, and which choice is made can have noticeable effects
  on what categorical constructions are permissible.  In this
  expository paper we summarize and compare a number of such
  ``set-theoretic foundations for category theory,'' and describe their
  implications for the everyday use of category theory.  We assume the
  reader has some basic knowledge of category theory, but little or no
  prior experience with formal logic or set theory.
\end{abstract}

\section{Introduction}
\label{sec:introduction}

Since its earliest days, category theory has had to deal with
set-theoretic questions.  This is because unlike in most fields of
mathematics outside of set theory, questions of \emph{size} play an
essential role in category theory.

A good example is Freyd's \emph{Special Adjoint Functor Theorem}: a
functor from a complete, locally small, and well-powered category with
a cogenerating set to a locally small category has a left adjoint if
and only if it preserves small limits.  This is always one of the
first results I quote when people ask me ``are there any real
\emph{theorems} in category theory?''  So it is all the more striking
that it involves in an \emph{essential way} notions like `locally
small', `small limits', and `cogenerating set' which refer explicitly
to the difference between sets and proper classes (or between small
sets and large sets).

Despite this, in my experience there is a certain amount of confusion
among users and students of category theory about its foundations, and
in particular about what constructions on large categories are or are
not possible.  Most introductory category theory books and courses
quite rightly ignore deeper set-theoretic questions, which will only
confuse most beginners.  However, intermediate and advanced students
of category theory may naturally begin to wonder about these
questions.

It turns out that there are several possible foundational choices for
category theory, and which choice is made can have noticeable effects
on what is possible and what is not.  The purpose of this informal
paper is to summarize and compare some of these proposed foundations,
including both `set-theoretical' and `category-theoretical' ones, and
describe their implications for the everyday use of category theory.
I assume the reader has some basic knowledge of category theory, such
as can be obtained from~\cite{maclane} or~\cite{awodey:cat-thy}, but
little or no experience with formal logic or set theory.  I found some
brief excursions into mathematical logic unavoidable, but I have tried
to explain all logical notions as they occur and relegate the more
complicated logical discussion to footnotes.

Related papers include~\cite[Appendix II]{feferman:fdns-of-ct} (by G.\
Kreisel) and~\cite{blass:ct-setth}; a historical survey can be found
in~\cite{kromer:tool-object}.  Since my goal is expository, I have
tried to restrict citations to those most useful for the student, when
such exist.  Likewise, this is not a history paper, so I will
concentrate on current understanding.

Finally, this is a work in progress, so suggested corrections and
improvements, both in mathematical content and exposition, are
extremely welcome.  Please also let me know if you have suggestions
for other references to include.  Numerous people have already helped
me with editing and feedback; I would like to especially thank Damir
Dzhafarov and Kenny Easwaran for useful discussions about set theory
and foundations; Antonio Montalban for explaining why the reflection
principle doesn't violate the incompleteness theorem; Peter Lumsdaine
for pointing out that a transitive set model without full power sets
can be countable; and Colin McLarty for directing my attention to the
categorical replacement axiom.

\section{Size does matter}
\label{sec:size-does-matter}

Before diving into set theory, it's natural to wonder why we need to
worry about size issues at all.  In this section I'll review two
theorems of basic category theory (interestingly, both due to Peter
Freyd) which, I think, display the essential nature of size
considerations.  Since this section is just motivation, I'll be vague
about the exact meaning of `set', `class', `small', and `large',
assuming the reader has some familiarity with their use.

First of all, we say that a category is \textbf{complete} if it admits
all limits indexed by \emph{small categories}; that is, categories
with only a set of objects and a set of morphisms.  The basic example
of a complete category is \Set: if \bA\ is small and $F\maps
\bA\to\Set$ is a functor, then the limit set $\lim(F)$ consists of
families $(x_a)_{a\in\bA}$ such that $x_a\in Fa$ and for all $f\maps
a\to b$ in \bA, $Ff(x_a)=x_b$.  There is the dual notion of
\textbf{cocomplete}.

It is essential in giving this definition that we restrict to
\emph{small} limits, since there are many large limits that \Set\ does
not admit.  For example, if $X$ is a set with more than one element,
then the $A$-fold product $\prod_{a\in A} X$ exists if $A$ is any set,
but not if $A$ is a proper class.  More generally, we have the
following, which is our first theorem in which size considerations are
essential.

\begin{thm}\label{thm:small-complete}
  If a category \bA\ has products indexed by the collection
  $\mathrm{Arr}(\bA)$ of arrows in \bA, then \bA\ is a preorder.  In
  particular, any small complete category is a preorder, and no large
  category that is not a preorder can admit products indexed by proper
  classes.
\end{thm}
\begin{proof}
  Suppose that we had two different arrows $f,g\maps a\to b$, and form
  the product $\prod_{\mathrm{Arr}(\bA)} b$.  Then $f$ and $g$ give us
  $2^{|\mathrm{Arr}(\bA)|}$ different arrows $a\to
  \prod_{\mathrm{Arr}(\bA)} b$; but there are only
  ${|\mathrm{Arr}(\bA)|}$ total arrows in \bA, a
  contradiction.\footnote{The use of proof by contradiction in this
    argument is essential.  In intuitionistic logic the theorem can
    fail; see~\cite{hyland:small-cplt}
    or~\cite[Ch.~24]{mclarty:ecat-etop}.}
\end{proof}

Thus, in order to capture most interesting examples, the notion of
\emph{complete category} must allow \emph{large} categories, but
restrict to \emph{small} limits.

However, many large categories \emph{do} admit \emph{some} large
limits.  For example, most large categories admit an intersection for
any family of monomorphisms with common codomain, no matter how large.
This is usually, but not always, because the category is
\textbf{well-powered}, meaning that each object has only a set of
isomorphism classes of subobjects.  Other large limits also often
exist; see, for example,~\cite{kelly:totality,kk:large-lim}.

Well-poweredness figures prominently in our second example: the
Special Adjoint Functor Theorem.  Recall that a family $Q$ of objects
in a category \bA\ is said to be \textbf{cogenerating} if whenever
$f\neq g\maps a\toto b$ are unequal parallel arrows, there exists an
arrow $h\maps b\to q$ with $q\in Q$ such that $hf\neq hg$.  We are
usually only interested in this when $Q$ is a set.  Recall also that a
category \bA\ is \textbf{locally small} if for any objects $a,b$ the
collection of morphisms $\bA(a,b)$ is a set.

\begin{thm}
  If \bA\ is locally small, complete, well-powered and has a
  cogenerating set, and \bB\ is locally small, then a functor $G\maps
  \bA\to\bB$ has a left adjoint if and only if it preserves small
  limits.
\end{thm}
\begin{proof}
  It suffices to construct, for each $b\in \bB$, an arrow $b\to GFb$,
  for some object $Fb\in \bA$, which is initial among arrows $b\to
  Ga$.  To define $Fb$, we first form the product $p=\prod_{q\in
    Q}\prod_{\bB(b,Gq)} q$ in \bA.  This product exists since $Q$ is a
  set, \bB\ is locally small, and \bA\ is complete.  Since $G$
  preserves products, we have an induced map $b\to Gp$.  Now let $Fb$
  be the intersection of all monomorphisms $a\mono p$ such that $b\to
  Gp$ factors through $Ga\to Gp$.  This intersection exists since \bA\
  is well-powered and complete.  Since $G$ preserves monomorphisms and
  intersections, we have an induced map $b\to GFb$.  We leave it to
  the reader to verify that this has the desired universal property
  (or see~\cite[V.8]{maclane}).
\end{proof}

For example, if \bA\ satisfies the hypotheses of the theorem and \bC\
is any small category, then the functor category $\pmb[\bC,\bA\pmb]$
is locally small and the diagonal functor $\Delta\maps
\bA\to\pmb[\bC,\bA\pmb]$ preserves limits, hence has a left adjoint.
Thus any such \bA\ is also \emph{cocomplete}.

Of course, as pointed out in the introduction, size distinctions play
an essential role in this theorem.  As stated, it applies to small
categories just as well as large ones, but it becomes somewhat
degenerate: any small category is locally small, well-powered, and has
a cogenerating set, so we obtain the following.

\begin{cor}
  If \bA\ is a complete lattice and $G\maps \bA\to\bB$ preserves
  greatest lower bounds, then it has a left adjoint.
\end{cor}

While undoubtedly important, this result is only a pale shadow of the
full Adjoint Functor Theorem.  Moreover, the Adjoint Functor Theorem
is not just a bit of fluff; there are examples even outside of pure
category theory where it is the only known way to construct an
adjoint.  So, like it or not, we are forced to deal with the question
of size in category theory.

\section{ZFC}
\label{sec:zfc}

With that motivation under our belts, we now turn to a quick summary
of set theory.  A natural question to begin with is: what is a
\emph{set}?  One modern answer is that \emph{sets} are special sorts
of \emph{collections}, which can be manipulated in well-defined ways
that
\begin{enumerate}[(a)]
\item suffice for applications in mathematics, but
\item are not powerful enough to reproduce the well-known paradoxes.
\end{enumerate}
There are three classical paradoxes of set theory, traditionally known
as Russell's, Cantor's, and Burali-Forti's.\footnote{Cesare
  Burali-Forti (1861--1931) is one of those mathematicians who are
  easily mistaken for two people by the unwary student.  Other
  distinguished members of this club include Tullio Levi-Civita
  (1873--1941) and G\"osta Mittag-Leffler (1846--1927).}  Russell's
paradox is non-categorical in flavor, while Burali-Forti's requires
ordinal numbers (see \autoref{thm:global-choice}), so here we recall
only Cantor's.

\begin{thm}
  There is no set containing all sets as members.
\end{thm}
\begin{proof}
  Suppose that $V$ were such a set.  Then every subset of $V$, being a
  set, would be a member of $V$; thus $\cP V \subset V$ and so
  $|\cP V| \le |V|$, contradicting Cantor's proof by diagonalization
  that $|A|<|\cP A|$ for any $A$.
\end{proof}

Thus, if we want to manipulate sets in the intuitive ways we are used
to, there must be some limitation on what collections we are allowed
to call `sets'.  The modern solution is to use a system of axioms
which allows us to construct enough sets to do mathematics, but not to
construct problematic sets such as $V$.

The set-theoretic axioms which have come to be accepted as standard
are today called \textbf{ZFC (Zermelo-Fraenkel set theory with
  Choice)} and can be found in any book on set theory (I
like~\cite{devlin:joy-of-sets,enderton:set-theory} as introductions,
while~\cite{jech:set-theory} is encyclopedic) or on the Internet.  For
later reference, we divide the axioms of \zfc\ into four types.
\begin{enumerate}
\item The \emph{basic} axioms: extensionality, foundation (or
  regularity), pairing, union, empty set, and separation.
\item The \emph{size-increasing} axioms: replacement and power set.
\item The \emph{size-assertion} axiom: infinity.
\item The axiom of \emph{choice}.
\end{enumerate}
Most of these are well-known or obvious.  The axiom schema of
\emph{separation} states that for any set $A$ and any definable
property $\ph(x)$, the set $\{x\in A \mid \ph(x)\}$ exists.  The axiom
schema of \emph{replacement} states that for any set $A$ and any
definable property $\ph(x,y)$ such that for any $x\in A$, there is a
unique $y$ with $\ph(x,y)$, the set $\{y\mid \exists x\in
A\;\ph(x,y)\}$ exists.  These are both both axiom \emph{schemas}:
there is one separation axiom for each definable property \ph\ and one
replacement axiom for each definable and `functional' property \ph.
In the presence of the other axioms, the replacement schema implies
the separation schema.

To be precise, \emph{definable} here means `definable in the formal
first-order language of set theory'.  But don't worry if you don't
know any logic; this really just means that it can be described in
ordinary mathematical language, referring only to sets or things that
can be defined in terms of sets (which includes most of mathematics).
For example, `$x$ is a continuous function from \bbR\ to \bbR' is a
definable property of $x$, so separation allows us to form the set of
all such functions (taking, $A$ to be, say, the power set of
$\bbR\times\bbR$).

All the ordinary constructions of mathematics can be performed using
these axioms.  For instance, we can define the ordered pair $(a,b)$ as
the set $\{\{a\}, \{a,b\}\}$, which exists by pairing, and
the cartesian product $A\times B$ as
\[A\times B = \Big\{ z \in \cP\big(\cP(A\cup B)\big) \;\Big|\;
\exists \,a\in A : \exists\, b\in B : z = (a,b)
\Big\},
\]
which exists by power-set and separation.  Similarly, the set $B^A$ of
functions from $A$ to $B$ can be defined by
\[B^A = \Big\{ f \in \cP(A\times B) \;\Big|\;
\forall\, a\in A: \exists\,!\, b\in B : (a,b)\in f
\Big\}.
\]

Let me comment briefly on the axiom schema of replacement, which may
seem the strangest one in the list from a categorical point of view.
In particular, it may seem odd to call it a size-increasing axiom,
since it merely replaces a set by an isomorphic one, or at most a
quotient.  We will see in later sections that given the other (also
non-categorical) axioms of \zfc, replacement in fact allows us to
construct much larger sets than would otherwise be possible.  But we
will also see that above and beyond this, replacement plays a subtle
and important role in category theory---so much so that this paper
could easily have been subtitled ``a tale of the replacement axiom''!

\begin{rmk}
  The approach of \zfc, and its relatives to be described in later
  sections, is not the only way to avoid paradoxes in set theory.  For
  example, in \emph{NFU (New Foundations with Urelements)}, any
  collection of things characterized by a `stratified' property is a
  set.  This allows for the existence of an actual set of all sets,
  while still avoiding paradoxes; see~\cite{holmes:nfu} for a good
  introduction.  However, \nfu\ is not much good for category theory,
  since the category \Set\ it produces is not cartesian
  closed~\cite{mclarty:nf-not-cc}.
\end{rmk}

\section{Ordinals and cardinals}
\label{sec:ordinals}

We now briefly review the theory of ordinal and cardinal numbers.
Succinctly, a \textbf{cardinal number} is a canonically chosen
representative for a bijection class of sets, while an \textbf{ordinal
  number} is a canonically chosen representative for an isomorphism
class of well-ordered sets.  Here a \emph{well-ordering} on a set is a
total ordering such that every nonempty subset has a least element.
The ordinals have the following properties.
\begin{enumerate}
\item Every ordinal \al\ has an immediate successor $\al+1$, obtained
  by adding an extra element at the end of a well-ordering of type \al.
\item There is a natural well-ordering on the collection of all
  ordinals: $\al\le\be$ iff \al\ is isomorphic to an initial segment
  of \be.
\item The induced well-ordering on $\{\be:\be<\al\}$ is in the
  isomorphism class represented by \al.\label{item:ord-neumann}
\item Every set is well-orderable, and hence bijective to some ordinal
  (this is equivalent to the axiom of choice).
\end{enumerate}
Because of \ref{item:ord-neumann}, one definition of an ordinal number
(due to von Neumann) is as the set of all smaller ordinals.  Thus, the
ordinals begin with the natural numbers $0=\emptyset$, $1=\{0\}$,
$2=\{0,1\}$, and so on, but continue afterwards with
\[\om, \om+1,\dots, \om\cdot 2, \dots, \om\cdot 3,\dots,
\om^2,\dots, \om^3,\dots, \om^\om,\dots, \om^{\om^\om},\dots,
\raisebox{.75cm}{\xymatrix@R=.5pc@C=.5pc{&\\\om^{\om^{\om}} \ar@{..}[ur]}},\dots
\]
We note in passing that the replacement axioms are first necessary to
construct $\om\cdot 2 = \{0,1,2,\dots, \om,\om+1,\om+2,\dots\}$.
Without replacement, we can construct each ordinal $\om+n$, but we
cannot collect them all as elements of a single set.  We can construct
well-ordered sets \emph{isomorphic} to $\om\cdot 2$ without
replacement, however, so the von Neumann definition of an ordinal is
only appropriate in the presence of replacement.

All the ordinals listed above are \emph{countable} (bijective with
\om).  We denote the first uncountable ordinal by $\om_1$, the next
ordinal not bijective with $\om_1$ by $\om_2$, and so on.  In fact, we
can define a cardinal number to be an ordinal not bijective with any
smaller ordinal.  It follows that the cardinal numbers are also
well-ordered, and can be indexed by the ordinal numbers.  We write
$\aleph_\al$ for the $\al^{th}$ cardinal number; thus $\aleph_0=\om$,
$\aleph_1=\om_1$, and so on.  Every set $X$ is bijective to a unique
cardinal $|X|$, called its cardinality.  If \ka\ is a cardinal, we
write $2^\ka$ for $|\cP\ka|$.

There are two types of ordinals.  Those of the form $\al+1$, for some
\al, are called \textbf{successor ordinals}, while the rest are called
\textbf{limit ordinals}.  Every cardinal is a limit ordinal.  A
cardinal is a \textbf{successor cardinal} or a \textbf{limit cardinal}
just when its indexing ordinal is a successor or limit ordinal.

Just as the well-ordering of \bbN\ justifies the usual sort of
mathematical induction, the well-ordering of ordinals justifies
definition and proof by \textbf{transfinite induction}.  This involves
proving or constructing something in stages, one for each ordinal.
The two cases of successor ordinals and limit ordinals are usually
dealt with differently; for $\al+1$ we base the construction on \al,
while for a limit ordinal \be\ we base it on all other ordinals
$\al<\be$.  We consider several important examples, three from set
theory and one from category theory.

\begin{eg}\label{thm:aleph-beth}
  The definition of the alephs can be phrased as a transfinite
  induction: we set $\aleph_0=\om$, let $\aleph_{\al+1}$ be the
  smallest cardinal greater than $\aleph_\al$, and for a limit $\be$
  we let $\aleph_\be = \lim_{\al<\be} \aleph_\al$.

  Similarly, we define $\beth_0=\om$, $\beth_{\al+1} = 2^{\beth_\al}$,
  and $\beth_\be = \lim_{\al<\be} \beth_\al$.  The \emph{Generalized
    Continuum Hypothesis (GCH)} is equivalent to
  $\aleph_\al=\beth_\al$ for all \al.
\end{eg}

\begin{eg}
  Define a set $V_\al$ for each ordinal \al\ by transfinite induction
  as follows.
  \begin{align*}
    V_0 &= \emptyset\\
    V_{\al+1} &= \cP V_\al\\
    V_\be &= \bigcup_{\al<\be} V_\al \qquad\text{($\be$ a limit)}.
  \end{align*}
  The sets $V_\al$ are called the \textbf{cumulative hierarchy}.  The
  axiom of foundation is equivalent to the assertion that every set is
  in $V_\al$ for some \al; this is often phrased as $V=\bigcup V_\al$,
  even though the class $V$ of all sets is not itself a set.

  The \textbf{rank} of a set $X$ is the smallest \al\ such that $X\in
  V_\al$.  For instance, the rank of each ordinal \al\ is $\al+1$.
  Most ordinary mathematical objects, as usually constructed from
  sets, have very low rank: the rank of $\bbN=\om$ is $\om+1$, the
  rank of \bbZ\ is $\om+2$, the rank of \bbQ\ is $\om+5$, and the rank
  of \bbR\ is $\om+9$ (or even less, if we are sufficiently clever).
\end{eg}

\begin{eg}\label{thm:v=l}
  For any set $X$, let $\Def(X)$ denote the set of all subsets of $X$
  which are \emph{definable} from elements of $X$.  By this I mean all
  sets of the form $\{x\in X \mid \ph(x)\}$ for some definable
  property $\ph(x)$ which refers only to elements of $X$---that is,
  its parameters and quantifiers (``for all $y$'' or ``there exists
  $y$'') range only over elements of $X$.  We define sets $L_\al$ by
  transfinite induction as follows.
  \begin{align*}
    L_0 &= \emptyset\\
    L_{\al+1} &= \Def(L_\al)\\
    L_\be &= \bigcup_{\al<\be} L_\al \qquad\text{($\be$ a limit)}.
  \end{align*}
  The sets $L_\al$ are called the \emph{constructible hierarchy} and
  the class $L=\bigcup L_\al$ is called the \emph{constructible
    universe}.  The \emph{axiom of constructibility} is the assertion
  that $V=L$; that is, that every set is constructible.  (Note that
  this does \emph{not} mean $V_\al=L_\al$!)  Most set theorists do not
  believe this axiom is `true', because it is so restrictive about
  what sets exist, but it cannot be proven or disproven from \zfc\
  alone, though it does contradict some large cardinal axioms (see
  \S\ref{sec:large-cardinals}).
\end{eg}

\begin{eg}\label{thm:soa}
  Let \bA\ be a cocomplete category and $S$ be a \emph{pointed
    endofunctor}, meaning a functor $S\maps \bA\to\bA$ equipped with
  a natural transformation $\si\maps \Id_\bA \to S$.  For any object
  $X\in \bA$ we define a sequence of objects $X_\al$, together with
  morphisms $X_\gm\to X_\al$ for $\gm<\al$, by transfinite induction
  as follows.
  \begin{align*}
    X_0 &= X\\
    X_{\al+1} &= S X_\al\\
    X_{\be} &= \colim_{\al<\be} X_\al \qquad\text{($\be$ a limit)}.
  \end{align*}
  The colimit at limit ordinals is, of course, over the diagram
  formed by the morphisms $X_\gm\to X_\al$, which we define by a
  parallel transfinite induction.  Namely, at a successor stage each
  morphism $X_\gm \to X_{\al+1}$ is the composite
  \[X_\gm \too X_\al \too[\si] SX_\al= X_{\al+1},
  \]
  while at a limit stage the morphisms $X_\gm \to X_{\be}$ are just
  the colimit cocone.

  If for some ordinal \de, the maps $X_\gm\to X_\al$ are isomorphisms
  for all $\de<\gm<\al$, we say that this process \emph{converges}.
  The intuition is that this happens when $SX$ is defined from $X$ by
  a `small amount of data', since then for a large enough limit
  ordinal \al, all the data necessary to construct $SX_\al$ will be
  contained in the objects $X_\gm$ for $\gm<\al$, so we will have
  $SX_\al\iso X_\al$.  Converging sequences of this sort are often
  used to construct reflections for subcategories and colimits in
  categories of algebras; an encyclopedic reference
  is~\cite{kelly:transfinite,kelly:transfinite-add}.

  A similar procedure is often followed in homotopy theory, but in
  this case usually we instead want the maps $X_\gm\to X_\al$ to
  become `weak equivalences' in an appropriate sense.  For example, if
  $\bA$ is the category of topological spaces and `weak equivalence'
  means `weak homotopy equivalence' (that is, a map inducing
  isomorphisms on all homotopy groups), then it usually suffices to
  take $\de=\om$.  This is because homotopy groups are detected by
  maps out of spheres, but spheres are compact, and so a map from a
  sphere into a well-behaved sequential colimit must factor through
  some finite stage.  However, in more complicated arguments, very
  large values of \de\ may be necessary.  In homotopy theory this is
  called the \emph{small object argument}, because it relies on the
  `smallness' of objects like spheres; see, for
  instance,~\cite{hovey:modelcats}.  For a version of the small object
  argument which does converge in the category-theoretic sense,
  see~\cite{garner:soa}.
\end{eg}

\section{Logic and Incompleteness}
\label{sec:logic}

A common mistake is to regard the axioms of \zfc\ as assertions only
about `the real' universe of sets, when in fact they are satisfied by
many different `universes of sets'.  This is not a philosophical
statement, but a mathematical one.  The reader is free to entertain a
Platonic belief that a `real' universe of sets exists (as many set
theorists seem to do), but it will still be true that the axioms of
\zfc\ support many different models in addition to this `real' one.  To
clarify the situation it is helpful to consider an analogy.

The \emph{axioms of group theory} are the following; they deal with a
collection of things and a binary operation ``$\cdot$''.
\begin{itemize}
\item For all $x,y,z$ we have $x\cdot (y\cdot z)=(x\cdot y)\cdot z$.
\item There exists an $e$ such that for all $x$, we have $x\cdot e = x
  = e\cdot x$.
\item For all $x$, there exists a $y$ such that $x\cdot y = e = y\cdot
  x$.
\end{itemize}
A \emph{model} of these axioms is a collection of things with a binary
operation satisfying them.  Of course, this is just a group.  We can
prove theorems from the axioms, which will then be true statements
about any group.  However, some statements, like ``for all $x$ and
$y$, $x\cdot y = y\cdot x$'', are neither provable nor disprovable
from the axioms; these are true for some groups and false for others.

In fact, \emph{G\"odel's Completeness Theorem} says that if a
statement is unprovable from the axioms of a theory \cT, then there
exist some models of the theory in which it is true and others in
which it is false.  We say that a theory is \textbf{consistent} if its
axioms do not imply a contradiction; the completeness theorem can then
be rephrased as ``any consistent theory has a model''.  Conversely,
the rather more obvious \emph{Soundness Theorem} says that any theory
with a model is consistent.

Now, the axioms of \zfc, which we summarized in \S\ref{sec:zfc}, deal
with a collection of things and a binary relation $\in$.  A
\emph{model} of the axioms of set theory is a collection of things,
which we usually call \emph{sets}\footnote{As we have stated them, the
  axioms of \zfc\ do not allow objects which are not sets.  It is easy
  to modify them to allow such `urelements', but there seems little
  point to doing so, since experience shows that everything in
  mathematics can be constructed using sets.}, together with a binary
relation $\in$, usually called \emph{membership}, which satisfy the
axioms.  Let us call such a model a \textbf{universe}.

We can prove many theorems from \zfc\ (in fact, we can develop most of
mathematics), and these theorems will then be true statements about
any universe.  However, just as for the theory of groups, some
statements are neither provable nor disprovable from the axioms; a
classic example is the Continuum Hypothesis (CH).  In fact, given any
universe, we can construct from it both a universe in which \ch\ is
true and a universe in which \ch\ is false.  The former is easy to
describe: the constructible universe $L=\bigcup L_\al$ is always a
model of \zfc+\ch.  The latter requires a more involved technique
called \emph{forcing} which is irrelevant to us here.  Thus, \emph{if}
\zfc\ is consistent, then both \zfc+\ch\ and \zfc\ + not \ch\ are
consistent.

Of course, one is entitled to wonder whether \zfc\ \emph{is} consistent;
that is, whether there are any universes.  This is no more or less
valid, from a purely logical standpoint, as wondering whether there
are any groups.  We are used to the existence of lots of groups, but
all the groups we are familiar with are constructed within the
framework of a stronger theory---namely, set theory.  In other words,
\emph{assuming} the existence of a universe, we can construct groups,
but with only the axioms of group theory, we can't expect to get
anywhere.  By analogy, we can't expect to be able to prove the
existence of a universe unless we work within the framework of some
yet stronger theory.

\emph{G\"odel's Second Incompleteness Theorem} is a formal way of
saying this: no reasonable\footnote{By `reasonable' I mean that there
  is a systematic way to verify whether or not any given statement is
  an axiom.  This excludes, for example, the system whose axioms are
  `all true statements about the natural numbers', to which of course
  the incompleteness theorem does not apply, but which is not much use
  as an axiom system in practice.} and consistent axiom system $\cT$
which includes arithmetic (such as \zfc) can prove its own
consistency.  It can be found in many books on logic (I
like~\cite{cl:logic} for mathematicians, while
everyone should read~\cite{geb}), but the proof is so simple in
outline that every mathematician should be exposed to it.  First, by
coding logical statements and proofs as natural numbers, G\"odel
enabled $\cT$ to talk about provability and consistency.  He then
constructed a statement $G$ \emph{about natural numbers} which said in
effect ``this statement is not provable in $\cT$''.  Thus, if $G$ is
provable, it is false.  Hence, if \cT\ is consistent, it cannot prove
$G$, and thus $G$ is true.  So there is a statement which is true, but
not provable in $\cT$; this is the First Incompleteness Theorem.  Note
that the completeness theorem then implies that any reasonable and
consistent theory has more than one model.

Now, the same coding of statements and proofs produces a statement
\emph{about natural numbers} which expresses `internally' our
inability to derive a contradiction from the axioms of \cT; call
this statement Con(\cT).  By internalizing the proof of the First
Incompleteness Theorem, we can prove \emph{in $\cT$} that
Con($\cT$) implies that $G$ is not provable, and hence that $G$ is
true.  Since $\cT$ cannot prove that $G$ is true, it follows that
$\cT$ cannot prove Con(\cT); this is the Second Incompleteness
Theorem.\footnote{The attentive reader will notice that this implies
  that if \cT\ is consistent, then there exist models of \cT\ in
  which Con(\cT) is false!  To make sense of this, remember that
  Con(\cT) actually says something like ``there does not exist a
  natural number $n$ which codes for a proof of $0=1$ in \cT.''
  Since \cT\ is consistent, there is no such proof, and thus no
  `real' natural number can code one, but bizarre models of \cT\ can
  contain `nonstandard' natural numbers which satisfy the
  \emph{arithmetical} property which we interpret as coding for such a
  proof.}

By the completeness theorem, it follows that we cannot prove in \cT\
the existence of a model for \cT.  Moreover, if we can prove the
existence of a model for \cT\ in some other theory $\cT'$, then $\cT'$
implies Con(\cT), and therefore Con(\cT) does not imply Con($\cT'$);
otherwise it would imply Con(Con(\cT)), contradicting the
incompleteness theorem.\footnote{Care is needed, however, when dealing
  with theories like \zfc\ that have infinitely many axioms.  It is
  possible to have a consistent theory \cT\ in which one can define a
  set $M$ and prove, for each axiom \psi\ of \cT, that \psi\ is true
  in $M$.  Such a theorem-schema does not contradict the
  incompleteness theorem; that would require proving a single theorem
  in \cT\ to the effect of ``for each axiom \psi\ of \cT, \psi\ is
  true in $M$''.  Ironically, one example of such a theory \cT\ is
  \zfc\ + ``\zfc\ is inconsistent'', which (by the incompleteness
  theorem) is consistent if \zfc\ is;
  see~\cite[IV.10]{kunen:set-theory}.  We will see in
  \S\ref{sec:nmod-refl} that \zfc\ itself is almost such a
  theory.\label{fn:inf-ax}} So if one theory can prove the existence
of a model for another, the first theory is irreducibly stronger in
this precise sense.

\begin{rmk}
  Actually, even without knowing the incompleteness theorem, or if we
  have an `unreasonable' \cT\ to which it does not apply, it is easy
  to see that a proof of Con(\cT) in \cT\ would be useless anyway.
  For since anything follows from a contradiction, if \cT\ were
  inconsistent, it would also prove Con(\cT).  Thus, a proof of
  Con(\cT) in \cT\ would still not allow us to conclude that \cT\ is
  \emph{actually} consistent.  The incompleteness theorem gives the
  stronger result that a proof of Con(\cT) in \cT\ in fact implies
  that \cT\ is \emph{in}consistent.
\end{rmk}

A similar result called \emph{Tarski's undefinability theorem} says
that there is no definable property \ph\ of natural numbers such that
$\ph(n)$ is true in some model of \cT\ if and only if $n$ is the
G\"odel code of a statement which is true in that model.  For suppose
there were.  Let $\psi_1,\psi_2,\dots$ be an enumeration of all
definable properties of natural numbers, and for any $n$ let $\#n(n)$
be the G\"odel code of the statement $\psi_n(n)$.  Now the statement
``$\ph(\#n(n))$ is false'' is a definable property of $n$, so it is
equal to $\psi_k$ for some $k$.  Then we have (in our model)
\[\psi_k(k) \iff \ph(\#k(k)) \text{ is false} \iff \psi_k(k) \text{ is
  false}\] which is absurd; thus \ph\ cannot exist.

On the other hand, since mathematical logic can be formalized in set
theory just as most branches of mathematics can, \zfc\ has no
problem talking about truth in a model which is a set.  That is, in
\zfc\ there is a definable property $\ph$ such that $\ph(n,x)$ is true
if and only if $n$ is the G\"odel code of a statement which is true in
$x$, regarded as a model of set theory.  This will also be important
later on.

\begin{rmk}
  The axioms of set theory and of group theory do differ in an
  important philosophical way.  The axioms of group theory are chosen
  because we see many objects `in nature' which satisfy them (whatever
  that means), and we want to study all these objects under one
  heading.  On the other hand, we do not see many examples of
  universes in nature---many people would argue that we see only one.
  The axioms of set theory are chosen for their usefulness,
  sufficiency, and consistency in working with sets, rather than
  claiming to be the `correct' description of an independently
  occurring class of models.
\end{rmk}

\section{Classes and large categories}
\label{sec:classes}

With some basic set theory under our belts, we now move on to category
theory.  By analogy with the theory of groups and the theory of sets
discussed in \S\ref{sec:logic}, we can consider the \emph{theory of
  categories}.  This theory deals with two types of things, called
`objects' and `arrows', together with domain, codomain, identity, and
composition functions satisfying unit and associativity axioms.  A
model of this theory is, of course, a category.  Note that this
abstract notion of `a category' can be defined without reference to
any sort of set theory.

In the context of a universe $V$, we generally refer to a category
whose collections of objects and arrows are sets in $V$ as a
\textbf{small} category.  When working with small categories with
respect to some universe, we have all the tools of set theory at our
disposal, and everything we might expect to be true, is.  For example,
given any two small categories $\bA$ and $\bB$, there is a small
category $\pmb[\bA,\bB\pmb]$ whose objects are functors $\bA\to \bB$
and whose arrows are natural transformations.

However, frequently even when working in the context of set theory, we
want to consider categories which are not small.  The obvious example
is the category $\Set = \Set[V]$ whose objects are all sets (that is,
all elements of the universe $V$) and whose arrows are all functions
between sets.  Cantor's paradox ensures that \Set\ is not small.
Non-small categories are usually called (surprise!) \textbf{large}.

If we are content to work with one large category at a time, we can
just use the theory of categories described above.  It is when we want
to construct new large categories and functors that we run into
problems, because the powerful tools of set theory are no longer at
our disposal for working with collections of objects that are not
sets.  Our goal is to consider various methods for dealing with this
problem.

First of all, there is an approach that remains completely within \zfc.
We define a \textbf{class} to be a collection of sets specified by
some property expressible in the language of set theory.  Instead of
working directly with classes, for which we have no axioms, we can
then work instead with the properties which characterize them.  For
example, the property ``$X$ is a pair $(G,\cdot)$ where $G$ is a set
and $\cdot$ is a binary operation on $G$ making it into a group'' is
expressible in the language of set theory, so there is a class of all
groups---by which we mean all groups defined from the universe $V$.
Of course, there is also a class of all sets, which we usually
identify with the universe $V$.  Note that classes of this sort are
actually implicit in the axioms of \zfc; for example, the axiom of
replacement says that the image of a set under any `class function' is
a set.  Some classes, of course, are sets; those which are not sets
are called \textbf{proper classes}.

We then define a \textbf{large category} to be one, such as \Set\ and
$\mathbf{Grp}$, whose collections of objects and arrows are classes in
this sense.  We can now perform some basic constructions on large
categories.  For example, if $P$ and $Q$ are properties expressible in
set theory, then ``$X$ is a pair $(Y,Z)$ such that $Y$ satisfies $P$
and $Z$ satisfies $Q$'' is also so expressible.  Thus the cartesian
product of two classes is a class, and the cartesian product of two
large categories is another large category.  We can also prove that
\Set\ and other familiar categories are complete and cocomplete, as
defined in \S\ref{sec:size-does-matter}.

\begin{rmk}
  Most large categories which arise in applications are also locally
  small; that is, they have only a set of morphisms between any two
  objects..  This property is undeniably important, as we saw in the
  proof of the Adjoint Functor Theorem, but we will mostly ignore it,
  because it plays almost no role when discussing foundations: locally
  small categories present exactly the same set-theoretic issues that
  all large categories do.
\end{rmk}

However, this approach to large categories has the disadvantage that
we have no axioms for manipulating classes; they are not `things' that
\zfc\ knows about at all.  Thus, instead of working with classes
directly, we have to work with the logical formulas which characterize
their elements, and interpret any construction in these terms.

In particular, the language of \zfc\ does not include a way to
\emph{quantify} over classes.  In other words, no theorem containing
the phrases ``for any large category \bA'' or ``there exists a large
category \bA'' can be even \emph{stated}, let alone proven, in \zfc.
This includes, for example, the Adjoint Functor Theorem.  Usually this
is dealt with by proving instead a `meta-theorem' of the form ``for
any large category \bA, we can prove in \zfc\ that (some statement about
\bA)'', but again in stating such a theorem we have moved beyond \zfc\
into some sort of `meta-language'.

Moreover, even this trick cannot handle theorems whose
\emph{hypotheses} involve quantification over classes.  For example,
consider the final statement in \autoref{thm:small-complete}: if a
large category has products indexed by proper classes, then it is a
preorder.  Even if we fix a large category \bA, the statement ``\bA\
has products indexed by proper classes'' is of the form ``for any
class \dots'', and thus cannot be stated in \zfc.  Hence we cannot
prove this result even as a theorem-schema in \zfc.

\section{Axioms for classes}
\label{sec:nbg}

To resolve these sorts of issues, we are motivated to extend \zfc\ by
introducing \emph{classes} as a new type of thing, in addition to
sets, along with axioms for manipulating them.  A good overview of
such \emph{class-set theories} can be found in~\cite{levy:classes}.
Probably the most common such theory is \textbf{von
  Neumann-Bernays-G\"odel (NBG) set theory}.  Its axioms can also be
found in books or on the Internet; as with \zfc\ we divide them into
several groups.
\begin{enumerate}
\item \emph{Typing}: only sets can be elements of sets or classes.
\item The \emph{basic} axioms: extensionality for sets and classes,
  the empty set, pairing and union for sets, foundation for sets and
  classes, and comprehension.
\item The \emph{size-increasing} axioms: power sets and limitation of
  size.
\item The \emph{size-assertion} axiom: an infinite set exists.
\item The axiom of \emph{choice}.
\end{enumerate}
The axiom schema of \emph{comprehension} states that for any property
$\ph(x)$ which does not quantify over classes, there is a \emph{class}
$\{x \mid \ph(x)\}$ of all \emph{sets} $x$ such that $\ph(x)$.  The
\emph{limitation of size} axiom says that a class is a set if and only
if it is not bijective with the class $V$ of all sets.  Thus, sets are
precisely the classes which are `not too big', while all proper
classes are the same size.

Comprehension and limitation of size easily imply separation and
replacement, so the sets in \nbg\ satisfy \zfc.  Moreover, \nbg\ can
be shown to be a \textbf{conservative} extension of \zfc.  That means
that any statement \emph{about sets} which is provable in \nbg\ is
also provable in \zfc.  In fact, if we start with any model of \zfc,
then taking the classes to be those defined in \S\ref{sec:classes}, we
obtain a model of \nbg.  Thus, using \nbg\ really entails no
`ontological' commitment beyond that of \zfc.

However, unlike \zfc, the language of \nbg\ allows us to quantify over
classes (although such quantifications cannot be used in the
comprehension axiom).  Thus, theorems such as the Adjoint Functor
Theorem can be stated and proven formally within \nbg.  To prove the
large-category version of \autoref{thm:small-complete}, we have to be
careful to avoid talking about $2^{|\mathrm{Arr}(\bA)|}$, which
doesn't exist in \nbg, but we can instead directly derive a
contradiction to the axiom of limitation of size.

Moreover, \nbg\ makes constructions on classes easier to deal with.
For example, comprehension proves easily that any two classes have a
cartesian product, and thus so do any two large categories.  We can
also perform more complicated constructions as long as they don't
produce things that are `too big'.

\begin{eg}
  Let \bA\ be a large category; we construct its
  \emph{idempotent-splitting} \bAbar, also called \emph{Karoubian
    completion} or \emph{Cauchy completion}.  The objects of \bAbar\
  are the idempotents of \bA; that is, arrows $e$ with $ee=e$.  The
  property ``$e$ is an arrow of \bA\ and $ee=e$'' is expressible in
  \nbg\ and doesn't quantify over classes, so by comprehension, there is
  a class of all such idempotents.  Similarly, the arrows of \bAbar\
  from $e$ to $e'$ are the arrows $f$ with $fe=f$ and $e'f=f$; this is
  likewise expressible without quantification over classes, so there
  is a class of all such arrows.  The rest of the structure of \bAbar\
  follows in the same way.
\end{eg}

\begin{eg}
  Let $\bA$ be a large monoidal category; we construct a strict
  monoidal category $\bA'$ monoidally equivalent to $\bA$.  One of the
  usual constructions is to let the objects of $\bA'$ be finite
  strings of objects of $\bA$, with morphisms induced by those of
  $\bA$.  Since the property ``$X$ is a function from some natural
  number $n$ to the class of objects of $\bA$'' does not quantify over
  classes, the class of all such functions exists; we take this to be
  the class of objects of $\bA'$.  Any such function (as a class of
  ordered pairs) is bijective with $n$ and thus is a set.  The rest of
  the structure is similar.
\end{eg}

Another example of the usefulness of having a good axiomatic system
for classes lies in our ability to make a large number of choices.
Suppose that $\bA$ is a category with finite products; it is usual to
make a choice of a product $X\times Y$ for each pair of objects
$X,Y\in \bA$ in order to define a product functor
$\bA\times\bA\to\bA$.  In particular cases there is usually a standard
choice of $X\times Y$, but to do this in general one needs an axiom of
choice for the objects of $\bA$.  If \bA\ is a large category, then
one needs an axiom of choice for \emph{classes}.  We call this the
\textbf{axiom of global choice}; it turns out to have the following
equivalent forms.
\begin{enumerate}[(i)]
\item We can choose an element from each of any class of
  nonempty sets.\label{item:gc1}
\item We can choose an element from each of any collection of nonempty
  classes.\label{item:gc2}\footnote{Note, however, that there is a
    standard trick in \zfc\ which enables us to choose an element from
    each of any \emph{set} of nonempty \emph{classes}.  To be precise,
    if we have a formula $\ph(x,y)$ such that for any $x\in X$ there
    is a $y$ with $\ph(x,y)$, we can define a function $f$ on $X$ such
    that $\ph(x,f(x))$ for all $x\in X$.  We do this by first
    considering, for each $x\in X$, the class of all sets $y$ of
    \emph{least rank} such that $\ph(x,y)$; this is a set since it is
    a subset of some $V_\al$.  We then apply the ordinary axiom of
    choice.}
\item The class $V$ of all sets can be well-ordered.\label{item:gc3}
\end{enumerate}
Global choice (and hence ordinary choice as well) is a consequence of
the axioms of \nbg, by the following observation of von Neumann.

\begin{thm}\label{thm:global-choice}
  In \nbg, there is a well-ordering of $V$.
\end{thm}
\begin{proof}
  The class $\Omega$ of ordinals is well-ordered.  Thus, if it
  were a set, it would itself be an ordinal; but then it would have a
  successor, which is absurd.  Thus $\Omega$ is not a set (this
  is Burali-Forti's paradox).  By limitation of size, $\Omega$
  is bijective to $V$, and thus $V$ acquires a well-ordering.
\end{proof}

If one objects to global choice despite the pleasing cleverness of
this argument, it is not hard to modify the axioms of \nbg\ so that they
no longer imply it.  Here we should also
mention~\cite{makkai:avoiding-choice}.  However, global choice is
implicitly used by many familiar categorical constructions.

\begin{eg}
  We have already noted that if \bA\ has products, then applying
  version~(\ref{item:gc2}) of global choice, we can choose a product
  $X\times Y$ for each pair $X,Y$ and thereby define a functor
  $\bA\times\bA\to\bA$.  Similar remarks, of course, apply to other
  limits and colimits, and other objects defined by universal
  properties, such as tensor products.
\end{eg}

\begin{eg}
  If $F\maps \bA\to\bB$ is a functor which is full, faithful, and
  essentially surjective, then by choosing for each $b\in \bB$ an
  object $Gb\in\bA$ and an isomorphism $FGa\iso b$, we can construct
  an inverse equivalence $G\maps \bB\to\bA$.
\end{eg}

\begin{eg}
  By choosing one object in each isomorphism class, we can show that
  any large category has a skeleton.
\end{eg}

\begin{eg}
  Let $W$ be a class of morphisms in a large category \bA; we want to
  construct a `localization' $\bA[W\inv]$ by formally adding inverses
  to the morphisms in $W$.  The objects of $\bA[W\inv]$ are the same
  as those of \bA, while its morphisms are supposed to be equivalence
  classes of zigzags
  \[\cdot \oot[\sim] \cdot \too\cdot \cdots \cdot \too \cdot \oot[\sim]\]
  of morphisms in \bA, where the backwards arrows are in $W$.
  Comprehension guarantees there is a class of such zigzags.  However,
  we cannot define the quotient of a class by an equivalence relation
  whose equivalence classes are proper classes---at least not in the
  usual way, since no class can be an element of another class.  But
  what we can do instead is use global choice to choose one zigzag
  from each equivalence class, thereby obtaining a class of morphisms
  for $\bA[W\inv]$.

  Of course, in general, $\bA[W\inv]$ defined in this way need not be
  locally small or at all amenable to computation.  In practice, there
  are usually alternate ways to construct $\bA[W\inv]$, which also
  show that in those cases it is locally small.
\end{eg}

\begin{eg}
  The \emph{left derived functors} of a functor $F$ evaluated at an
  object $X$ of some abelian category are given by choosing a
  projective resolution $\dots\to P_1\to P_0 \to X$ and computing the
  homology of the chain complex $\dots\to FP_1\to FP_0$.  Of course,
  defining the whole derived functor requires the global choice of a
  projective resolution for each object $X$.
\end{eg}

On the other hand, \nbg\ is not quite as comfortable an axiom system as
we might like.  Consider mathematical induction, which is surely a
basic notion in mathematics if ever there was one.  In \zfc, we can
prove induction for any definable statement $\ph(n)$, by considering
the set $\{n \in\bbN \mid \text{not } \ph(n)\}$ and using the fact
that \bbN\ is well-ordered.  In \nbg\ the same argument works only if
$\ph$ does not involve quantification over class variables, due to the
analogous restriction in the comprehension axiom.  This includes all
statements for which \zfc\ could prove induction (as it must, since \nbg\
is a conservative extension of \zfc), but not all statements in our
`new language' which can refer to classes.

One might object to this consequence of \nbg\ on the philosophical
grounds that mathematical induction `should' be true for \emph{all}
statements $\ph$, without needing technical restrictions on
quantification.  But the failure of induction also has real
consequences for dealing with classes.  For instance, we cannot prove
by induction on $n$ that every large category \bA\ has an $n$-fold
cartesian product $\bA^n$.  In this case, we can construct $\bA^n$
directly as the class of functions from $n$ into \bA, but it is
troubling that the induction proof is not allowed in \nbg.\footnote{In
  defense of \nbg, I should say that it was originally conceived not
  to deal with large categories, but to provide a finitely
  axiomatizable theory equivalent to \zfc, and at that it succeeds.
  This is not immediately obvious, since we have stated comprehension
  as an axiom schema, but it turns out that a finite number of its
  instances suffice to imply the rest.  There is no great mystery
  about this: we simply observe that any definable property \ph\ is
  built up from a finite number of building blocks like `and', `or'
  and `there exists', and the corresponding class $\{x\mid \ph(x)\}$
  can be built up by a corresponding finite number of constructions
  like intersection, union, and projection.  However, it does depend
  on limiting comprehension to properties not quantifying over
  classes; thus \mk\ (see below) is not finitely axiomatizable in this
  way.}  If any reader can think of a natural statement about large
categories, normally proven by induction, and not having an obvious
alternate proof in \nbg, I would be very interested.

We may, of course, strengthen the comprehension axiom of \nbg\ to
allow formulas with arbitrary quantification (sometimes called
\emph{impredicative comprehension}), thereby recovering full
mathematical induction.  This gives a variant of \nbg\ usually called
\textbf{Morse-Kelley (MK) set theory}.  However, by doing so we lose
conservativity over \zfc.  In \mk\ one can prove that the class $V$ of
all sets is a model for \zfc\ and conclude that \zfc\ is consistent.
For example, to prove the separation and replacement axioms, we first
apply comprehension to construct the desired set or function as a
class, then apply limitation of size to conclude that it is a set.  It
then follows from the incompleteness theorem that Con(\zfc) does not
imply Con(\mk).\footnote{Why, the reader may reasonably wonder, can we
  not do the same in \nbg?  We can prove in \nbg\ that $V$ satisfies
  any \emph{particular} axiom of \zfc, but as in
  footnote~\ref{fn:inf-ax}, to conclude Con(\zfc) we need instead the
  single theorem ``for all axioms \psi\ of \zfc, $V$ satisfies \psi''.
  To even state this formally when there are infinitely many axioms,
  we must encode axioms by their G\"odel numbers.  Just as \zfc\ can
  talk about truth in set models, in \nbg\ we have a definable
  property \ph\ such that $\ph(n)$ is true if and only if $n$ is the
  G\"odel code of a true statement involving only sets---but this \ph\
  involves quantification over classes.  Thus, we require the strength
  of \mk\ to form the class $\{x \mid \psi(x)\}$ given only the
  G\"odel number of \psi, as is necessary for the above proof of
  Con(\zfc).  This distinction is well explained
  in~\cite{mostowski:impred,mostowski:cor-impred}, along with
  resulting concrete examples of the failure of full mathematical
  induction and full class comprehension in \nbg.}  So unlike \nbg,
\mk\ is a genuinely stronger theory than \zfc.

Even \mk, however, is not fully satisfactory in the constructions it
allows for large categories.  For example, if \bA\ and \bB\ are large
categories, nothing we have seen so far allows us to construct a
functor category $\pmb[\bA,\bB\pmb]$.  There is no problem when \bA\
is small, since then each functor $\bA\to\bB$ is itself a set by
replacement, and so we have a class of such---but when \bA\ is large,
this argument fails.  However, there is no intuitive reason preventing
us from making such a construction: the collection of functions from
one class to another seems like a perfectly good collection.

We could envision adding more axioms which enable us to perform these
and other constructions with classes.  In fact, the best possible
world would be if classes could be manipulated \emph{just like} sets,
and any construction we could do for sets could also be done for
classes.  This would be easy to achieve: we could just write down
another copy of the \zfc\ axioms, substituting `class' for `set'
everywhere in the second copy.  On the other hand, it seems terribly
wasteful to have two copies of every axiom, when all we really want to
say is that classes and sets behave in just the same ways, except that
sets can't be too large.  In the next section we consider a cleaner
solution.

\section{Inaccessibles and Grothendieck universes}
\label{sec:universes}

In \zfc, there are three ways to prove the existence of larger and
larger sets.
\begin{enumerate}[(1)]
\item By \emph{fiat}: the axiom of infinity asserts that there exists
  an infinite set.  Without it, only finite sets can be constructed.
\item By \emph{powers}: the axiom of power sets produces a set $\cP
  A$ larger than a given set $A$ (by Cantor's diagonalization
  argument).
\item By \emph{limits}: The axiom of replacement guarantees that the
  union of any family of sets indexed by a set is also a
  set.\footnote{In more detail, let $A$ be a set and $F$ a class
    function on $A$.  Such an $F$ is characterized by some definable
    property $\ph$ such that for any $x\in A$ there exists a unique
    $y$ with $\ph(x,y)$.  Then replacement gives the set $\{y \mid
    \exists x\in A\; \ph(x,y)\}$, to which we can then apply the union
    axiom to give the set $\bigcup_{x\in A} F(x)$.  It is natural to
    wonder why I call replacement the culprit here, when the union
    axiom seems at least equally culpable; one answer is that, as we
    will see below, small models of \zfc\ satisfying the union axiom
    abound, while ones satisfying replacement are quite rare.}  If the
  family is infinite and increasing in size, its union will be larger
  than any of its elements.  For example, this applies to a family
  such as $A$, $\cP A$, $\cP\cP A$, \dots.
\end{enumerate}
This is where I got the terminology `size-assertion' and
`size-increasing' in \S\ref{sec:zfc}; the axiom of infinity produces a
large set by fiat, while the axioms of power set and replacement
produce larger sets from existing ones.  However, not all
cardinalities can be reached by these methods; we introduce special
names for those that can't.
\begin{enumerate}
\item A cardinal is \textbf{uncountable} if it is larger than the
  smallest infinite set.
\item A cardinal \ka\ is a \textbf{strong limit} if for any $\la<\ka$ we
  have $2^\la < \ka$.
\item A cardinal \ka\ is \textbf{regular} if it is not the union of a
  family of sets of size $<\ka$ indexed by a set of size $<\ka$.
\item A cardinal is \textbf{inaccessible} if it is uncountable, a
  strong limit, and regular.
\end{enumerate}

For example, the first uncountable cardinal $\aleph_1$ is regular,
since the countable union of countable sets is countable, but it is
not a strong limit, since $2^{\aleph_0} \ge \aleph_1$.  On the other
hand, the cardinal $\beth_\om$ (see \autoref{thm:aleph-beth}) is a
strong limit, but not regular.


Now, in any universe $V$, the set $V_\al$ with its induced relation of
membership is itself a model of many of the axioms of \zfc.  It is
easy to see that if \al\ is a limit ordinal greater than $\om$, then
$V_\al$ satisfies the basic axioms, choice, power set, and
infinity---all the axioms of \zfc\ except replacement.  In particular,
since $|V_{\om+\al}| = \beth_{\al}$, without replacement we can only
construct sets of cardinality $<\beth_\om$.

It turns out that if \al\ is an inaccessible cardinal, then $V_\al$ is
also a model of replacement, and hence of all of \zfc.  (The converse,
however, is false, as we will see in \S\ref{sec:nmod-refl}.)  The
proof is the same as the proof in \mk\ that $V$ satisfies replacement,
using the inaccessibility of \al\ in place of limitation of size.
When \al\ is inaccessible we call $V_\al$ a \textbf{Grothendieck
  universe}.  One can equivalently define a Grothendieck universe to
be a set $U$ which is transitive ($x\in y\in U$ implies $x\in U$) and
closed under pairing, power sets, and indexed unions.  It turns out
that this is equivalent to asserting $U=V_\ka$ for some inaccessible
\ka; see~\cite{bourbaki:universe}.

Now suppose there exists an inaccessible, and let \ka\ be the smallest
inaccessible.  Then $V_\ka$ satisfies \zfc\ and also ``there does not
exist an inaccessible''; hence it is impossible to prove in \zfc\ that
there exists an inaccessible.  But this is not really surprising,
since we have essentially \emph{defined} inaccessibles to be those
cardinals which are unreachable by all the ways that \zfc\ knows to
build bigger sets!

More that this is true, however.  If we write \textbf{ZFC+I} for \zfc\
+ ``there exists an inaccessible'', then even assuming that \zfc\ is
consistent, it is not possible to prove that \zfci\ is consistent.
For just as we can prove in \mk\ that $V$ is a model of \zfc, we can
prove in \zfc\ that any Grothendieck universe is a model of \zfc;
hence \zfci\ implies Con(\zfc).  Thus, by the incompleteness theorem,
Con(\zfc) does not imply Con(\zfci).  In other words, in contrast to
the situation for \ch, we cannot construct, from an arbitrary
universe, another universe satisfying \zfci.  Note that Con(\zfc),
which is provable in \zfci\ but not in \zfc, is a statement not about
sets but about \emph{natural numbers}---albeit a rather complicated
one.

Laying aside questions of existence and consistency for the moment, we
can solve the problems raised at the end of \S\ref{sec:nbg} as
follows.  Working in \zfci, we choose an inaccessible \ka, and
re-define \emph{set} to mean `element of $V_\ka$' and \emph{class} to
mean `set not necessarily in $V_\ka$'.  Thus defined, sets and classes
will behave in exactly the same way, except that sets are limited in
rank.  A more common terminology, however, which we will adopt, is not
to redefine `set' but to refer to elements of $V_\ka$ as \textbf{small
  sets} and other sets as \textbf{large sets}.

Note that small sets have a limitation on \emph{rank} rather than
\emph{cardinality}.  No small set can be larger in cardinality than
\ka, but many sets with small cardinality are not small, such as the
singleton $\{\ka\}$.  The class of sets with small cardinality is not
a model of \zfc, since it is not closed under unions ($\bigcup \{\ka\}
= \ka$), nor is it itself a set (even a large one).  Thus, if we want
a category \Set\ of small sets whose collection of objects is a large
set, we need to use rank rather than cardinality.\footnote{There is,
  however, an approach to universes based on cardinality: for any
  infinite cardinal \ka\ let $H_\ka$ denote the set of sets which are
  \emph{hereditarily} of cardinality $<\ka$; that is, their transitive
  closure has cardinality $<\ka$.  Then $H_\ka\subset V_\ka$, and if
  \ka\ is regular and uncountable, $H_\ka$ satisfies all the axioms of
  \zfc\ except power set.  Moreover, \ka\ is inaccessible if and only
  if $H_\ka$ satisfies all of \zfc, and if and only if $H_\ka=V_\ka$.}

Of course, any category of sets is determined up to equivalence by the
cardinalities of its objects, so in some sense, this doesn't matter.
For this reason, it is common to also call sets \emph{small} if their
cardinality is small, and in category theory, where we work with
objects up to isomorphism anyway, this is usually harmless.

This third approach allows many more large categories than the first
two.  For instance, for any two large categories \bA\ and \bB, there
is a functor category $\pmb[\bA,\bB\pmb]$.  Note that unlike in \nbg,
where all proper classes have the same size, $\pmb[\bA,\bB\pmb]$ will
generally be larger than \bA\ or \bB.  In particular, the term
\emph{large category} in this approach is significantly more inclusive
than it is in the previous two.

One should think of the classes in \nbg\ and \mk\ as corresponding to
the large sets of rank (or cardinality) $\le \ka$.  In fact, one can
make this precise.  If \ka\ is inaccessible, we obtain a model of \mk\
(and hence \nbg) by taking $V_\ka$ for the sets and $V_{\ka+1} = \cP
V_\ka$ for the classes.  (Note that this means we can prove Con(\mk)
in \zfci, so \zfci\ is strictly stronger than \mk.)  Similarly, we
obtain a model of \nbg\ (but not \mk) if we take $V_\ka$ for the sets
and $\Def(V_{\ka})$ for the classes (recall \autoref{thm:v=l}).

In~\cite{street:topos}, sets and categories of cardinality $\le \ka$
(which are those appearing, up to isomorphism, in $V_{\ka+1}$) were
called \textbf{moderate}.  We will refer to sets and categories in
$\Def(V_{\ka})$ as \textbf{small-definable}.  Small-definable
categories are those which would exist as classes even in \zfc\ or \nbg;
this includes nearly all large categories in which we are actually
interested for their own sake.  However, the existence of other large
categories is frequently quite useful for formal reasons; it is the
primary advantage of using inaccessibles.  In the following examples
we assume an inaccessible \ka.

\begin{eg}
  Any large category $\bA$ has a presheaf category
  $\pmb[{\bA\op},\Set\pmb]$, and if \bA\ is locally small, it has a
  Yoneda embedding $y\maps \bA\into \pmb[{\bA\op},\Set\pmb]$.  Limits
  and colimits in large functor categories can still be calculated
  pointwise, so $\pmb[{\bA\op},\Set\pmb]$ is still complete and
  cocomplete.  Important properties of \bA\ can be expressed in terms
  of $y$; for example, \bA\ is \emph{total}
  (see~\cite{kelly:totality}) if $y$ has a left adjoint.  Totality can
  be expressed without reference to $\pmb[{\bA\op},\Set\pmb]$
  (although not without quantifying over classes), but at the price of
  a certain amount of economy and clarity.
\end{eg}

\begin{eg}
  Any large \bA\ has an endofunctor category $\pmb[\bA,\bA\pmb]$,
  which is strict monoidal under composition of functors.  If \bA\ is
  also monoidal, then the functor $\bA\to \pmb[\bA,\bA\pmb]$ taking
  $X\in\bA$ to $(X\ten -)$ is strong monoidal.  We can thus apply the
  idea of~\cite[XI.3 Ex.~3]{maclane} to find a strictification of \bA\
  inside $\pmb[\bA,\bA\pmb]$.
\end{eg}



\begin{eg}
  It is well-known in algebraic geometry that the category
  $\mathbf{Sch}$ of \emph{schemes} is equivalent to a certain
  subcategory of $\pmb[\mathbf{Ring},\Set\pmb]$.  The category
  $\mathbf{Sch}$ also has other definitions (for instance, as locally
  affine locally ringed spaces), showing that it is small-definable
  (up to equivalence) and locally small.  However, identifying it with
  a subcategory of $\pmb[\mathbf{Ring},\Set\pmb]$ is often useful, yet
  impossible unless the latter category exists.
\end{eg}

\begin{eg}
  Similarly, the category $\operatorname{Ind}(\bA)$ of
  \emph{ind-objects} in a large category \bA\ can be identified with
  the small filtered colimits of representable presheaves in
  $\pmb[\bA\op,\Set\pmb]$.  Like $\mathbf{Sch}$, the category
  $\operatorname{Ind}(\bA)$ has an alternate description showing that
  it is small-definable (up to equivalence, as long as \bA\ is), but
  this description of it is frequently also useful.
\end{eg}

\begin{eg}
  Let $U\maps \mathbf{CptHaus}\to\Set$ be the forgetful functor from
  compact Hausdorff spaces to sets.  A \emph{quasi-topological space}
  is a set $X$ equipped with a subfunctor of $\Set(U-,X)\maps
  \mathbf{CptHaus}\op\to\Set$ satisfying certain natural conditions;
  see~\cite{quasi-topologies}.  The category $\mathbf{QTop}$ of
  quasi-topological spaces is cartesian closed, and was a
  contender for a convenient category of spaces before the current
  ascendancy of compactly generated spaces (for which
  see~\cite[VII.8]{maclane} and~\cite[Ch.~5]{may:concise}).

  Since a single quasi-topological space contains a large amount of
  data, $\mathbf{QTop}$ is \emph{not} small-definable, though it is
  locally small.  Hence \nbg\ and \mk\ are insufficient to guarantee
  that $\mathbf{QTop}$ even exists.  Also, a fixed set $X$ supports a
  large number of quasi-topologies, and $\mathbf{QTop}$ is not
  well-powered or well-copowered.  It is, however, complete and
  cocomplete, and admits intersections of arbitrary families of
  monomorphisms and cointersections of arbitrary families of
  epimorphisms.
\end{eg}

\begin{eg}
  In addition to its reassuring psychological effect, moderateness
  (and small-definability) of a category can have mathematical
  consequences.  This is because a set $A$ is moderate if and only if
  we can express it as an increasing union of small sets indexed by
  small ordinals, $A = \bigcup_{\al<\ka} A_\al$ where each $A_\al$ is
  small.  Thus, we can prove things about moderate sets by a
  transfinite induction in which every stage of the induction is
  small.  For example, a proof of Freyd given in~\cite{street:topos}
  shows that if \bA\ is moderate, total, and the left adjoint of $y$
  preserves finite limits, then \bA\ actually has a \emph{small}
  generating set and is a Grothendieck topos.
\end{eg}

Sometimes it is useful to assume more than one inaccessible.  For
example, when doing formal category theory we may want to form the
category (or 2-category) of all large categories, or of all locally
small categories.  Of course, the class of all large categories is not
a set, even a large one.  (The class of small-definable categories is
a large set, but it is not closed under constructions such as functor
categories.)  To resolve this issue, we can use the techniques of
\S\S\ref{sec:classes}--\ref{sec:nbg} to introduce \emph{classes} that
are larger than large sets, or we can assume a second inaccessible
$\la>\ka$, define \textbf{large} to mean an element of $V_\la$, and
use \textbf{very large} to mean a set not necessarily in $V_\la$.
(Some authors have used `quasi-category' for what we call a `very
large category', but we eschew that term in view of its quite
different recent connotations~\cite{joyal:q_kan}.  The term
`meta-category' is also sometimes used.)

Having a very large category $\mathbf{CAT}$ of large categories allows
us to make statements like the following.
\begin{itemize}
\item Taking a small category to its presheaf category is a functor
  from the category $\mathbf{Cat}$ of small categories to the category
  $\mathbf{CAT}$ of large ones.
\item Taking a ring $R$ to the category of $R$-modules is a functor
  from $\mathbf{Ring}$ to $\mathbf{CAT}$.
\item Taking a monoidal category \bV\ to the category
  $\bV\text{-}\mathbf{Cat}$ of small \bV-enriched categories is a
  functor from $\mathbf{MONCAT}$ to $\mathbf{CAT}$.
\end{itemize}
If we want to have a functor taking \bV\ to the very large category
$\bV\text{-}\mathbf{CAT}$ of \emph{large} \bV-enriched categories, its
codomain will have to be an \emph{extremely large} category of very
large categories, so we need at least \emph{three} inaccessibles.  One
is unavoidably reminded of how, in the original category theory
paper~\cite{em:first-ct}, (large) categories were introduced only in
order to serve as the domains and codomains of functors.

\section{Aside: large cardinals}
\label{sec:large-cardinals}

Notice the strong analogy between the axiom ``there exists an
inaccessible'' and the axiom of infinity: both construct ``by fiat''
larger sets than can otherwise be shown to exist.  Other axioms of
this sort, asserting the existence of inaccessibles satisfying various
extra properties, are also used in modern set theory.  Of course,
these are stronger assumptions than the mere existence of an
inaccessible.  In fact, the existence of even one cardinal with one of
these stronger properties usually implies the existence of many
smaller inaccessibles.

Let me attempt to give a flavor of just how large such \emph{large
  cardinals} can get.  (The terminological collision between `large
category' and `large cardinal' is unfortunate, but context usually
suffices to distinguish.)  A subset $X\subset\ka$ is said to be
\textbf{closed unbounded} if $\sup(X) = \ka$ and whenever $Y\subset X$
and $\sup(Y) <\ka$, then also $\sup(Y) \in X$.  A subset $X\subset
\ka$ is \textbf{stationary} if it has nonempty intersection with every
closed unbounded set.  Evidently any stationary set is unbounded, and
hence (if \ka\ is inaccessible) has cardinality \ka.

An inaccessible \ka\ is said to be \textbf{Mahlo} if the set of
inaccessibles less than \ka\ is stationary in \ka.  This implies that
there are \ka\ inaccessibles less than \ka, but also that there are
\ka\ inaccessibles $\la<\ka$ such that there are \la\ inaccessibles
less than \la, and \ka\ inaccessibles below \ka\ with this property,
and so on.  This pattern continues for many larger types of cardinals:
\begin{itemize}
\item If \ka\ is \emph{weakly compact}, then Mahlo cardinals are
  stationary in \ka.
\item If \ka\ is \emph{measurable}, then weakly compact
  cardinals are stationary in \ka.
\item If \ka\ is \emph{supercompact}, then measurable cardinals are
  stationary in \ka.
\item If \ka\ is \emph{extendible}, then supercompact cardinals are
  stationary in \ka.
\item If \ka\ is \emph{superhuge}, then extendible cardinals are
  stationary in \ka.
\end{itemize}
Unlike for Mahlo cardinals, these stationarity properties are not the
definitions of these larger cardinals, but consequences
thereof.\footnote{When comparing large cardinal axioms in set theory,
  \emph{consistency strength} is usually more important than raw size.
  Obviously, if there are many $X$ cardinals below any $Y$ cardinal,
  then the consistency of a $Y$ implies that of an $X$, but the
  converse is not always true.  For example, a huge cardinal implies
  the consistency of extendible cardinals, and hence of supercompact
  ones, but the least huge cardinal is less than the least
  supercompact cardinal (if both exist).}  Their actual definitions
can be found in books (for
example,~\cite{jech:set-theory,kanamori:higher-infinite}) or the
Internet.

The purpose of large cardinals in set theory is not just to see how
large sets can get, but to provide a yardstick to measure the strength
of other axioms, and in some cases even to prove new results about
`small' sets.  For instance, the existence of a measurable cardinal
implies that $V\neq L$, the existence of infinitely many Woodin
cardinals implies \emph{projective determinacy}, and at least some set
theorists hope that a large cardinal axiom can be found which will
settle the continuum hypothesis;
see~\cite{believing-axioms-i,believing-axioms-ii}.

Perhaps surprisingly, set theorists currently believe that there is an
upper bound to how large large cardinals can get: there are notions of
\emph{$n$-huge} cardinal for all $n<\om$, but the limiting case of an
`\om-huge' cardinal is known to be inconsistent with \zfc.  However,
even large-cardinal axioms only slightly weaker than \om-hugeness have
so far resisted all efforts to disprove them.

Returning to category theory, it is natural to wonder whether large
cardinal properties could have noticeable effects on \Set\ and
categories constructed from it.  This is indeed the case, although for
most category-theoretic purposes one inaccessible is as good as
another.  For example, it is proven in~\cite[A.5]{ar:loc-pres} that
$\Set\op$ has a small dense full subcategory if and only if there do
not exist arbitrarily large measurable cardinals.  If there are
\emph{no} measurable cardinals, then the single set \bbN\ is dense in
$\Set\op$, and thus $\Set\op$ is equivalent to a full subcategory of
the category of $M$-sets, where $M$ is the monoid of endomorphisms of
\bbN.\footnote{I am inclined to regard this as an argument in favor of
  the existence of measurable cardinals, though of course others may
  disagree.  Note that $\Set\op$ is also equivalent to the category of
  profinite Boolean algebras (by Stone duality), and to the category
  of algebras for the double-power-set monad on \Set.}

In general, while large cardinal axioms in set theory usually assert
the \emph{existence} of one or more cardinals with a certain property,
what tends to matter for category theory is the character of the
\emph{particular} `size of the universe' cardinal \ka.  For example,
what matters for the result quoted above is whether the measurable
cardinals are unbounded below \ka, rather than how many measurable
cardinals there are in `absolute' terms.  Moreover, in most cases, the
assertion that ``the cardinal of the universe has property $P$'' can
be phrased in \nbg\ (and sometimes even \zfc), without requiring the
existence of any sets larger than the universe.

The most interesting examples of this sort that I know of concern
\textbf{Vop\v{e}nka's principle}.  This has many equivalent forms;
here are a few categorical ones.
\begin{itemize}
\item No locally presentable category has a large discrete full
  subcategory.
\item Every complete or cocomplete category with a small dense full
  subcategory is locally presentable.
\item Every category with a small dense full subcategory is
  well-copowered.
\end{itemize}
None of these can even be stated in \zfc, since they all involve
quantification over large categories, but there is no problem in \nbg.
{Vop\v{e}nka's principle} has many other pleasing consequences for the
structure of locally presentable categories
(see~\cite[Ch.~6]{ar:loc-pres}), and also implies the existence of
arbitrary cohomological localizations in homotopy
theory~\cite{css:large-cardinal}, which are not known to exist in \zfc.

We say that a cardinal \ka\ is \textbf{Vop\v{e}nka} if {Vop\v{e}nka's
  principle} holds in \zfci\ with \ka\ as the size of the universe.
This is equivalent to saying that {Vop\v{e}nka's principle} holds in
$V_\ka$ regarded as a model of \mk, but stronger than the analogous
assertion involving \nbg.  Since $\Set\op$ is not locally presentable,
the assertions I have made so far imply that measurable cardinals are
unbounded below any Vop\v{e}nka cardinal, but more is true: if \ka\ is
Vop\v{e}nka, then the sets
\begin{gather*}
  \{\la<\ka \mid \la \text{ is measurable}\}\\
  \{\la<\ka \mid \la \text{ is extendible in } V_\ka\}
\end{gather*}
are stationary in \ka.  This makes the existence of even one
Vop\v{e}nka cardinal quite a strong assumption, as large-cardinal
axioms go.  On the other hand, Vop\v{e}nka cardinals are stationary in
any `almost huge' cardinal.

\section{Inaccessibles or not?}
\label{sec:phil-inacc}

So, should we assume inaccessibles?  As we have seen, the existence of
an inaccessible---even the \emph{consistency} of the existence of an
inaccessible---is unprovable from \zfc, so such an assumption is a
genuine strengthening of the axioms.  On the other hand, there are
many philosophical and mathematical arguments that can be advanced in
its favor; see~\cite{believing-axioms-i} for a fascinating discussion.
Moreover, we have seen that it is quite weak compared to many
large-cardinal axioms commonly used in modern set theory, and from
which no contradiction has yet been derived.  Thus, it seems unlikely
that the existence of inaccessibles can be disproven from \zfc\
(though it is provably impossible to prove that it can't be!).

However, from the point of view of ordinary category theory, these
questions are not as relevant, because the role of inaccessibles in
category theory is quite different from their role in set theory.  In
category theory, inaccessibles mostly play the role of a
\emph{convenience} which simplifies the statements and proofs of our
theorems, without really entailing any deep ontological commitment.
This is because we think of the small sets as the `real' sets; we only
introduced the large ones as a well-behaved model for proper classes.
All ordinary mathematical objects, like groups, rings, topological
spaces, manifolds, and so on, are small.

Moreover, all categories of ordinary objects which arise in practice,
such as \Set\ and $\mathbf{Grp}$, are small-definable, and would exist
as classes even in \zfc\ or \nbg.  And we saw that categories such as
$\mathbf{Sch}$ and $\operatorname{Ind}(\bA)$ have equivalent forms
which are small-definable, even though they can also be usefully
characterized with reference to categories that are not.

In a similar way, a statement like
``$\bV\mapsto\bV\text{-}\mathbf{Cat}$ is a functor'' can be regarded
merely as a `code' which encapsulates many individual statements in a
concise way.  For example, it implies that any monoidal functor
$\bV\to\bW$ induces a functor $\bV\text{-}\mathbf{Cat}\to
\bW\text{-}\mathbf{Cat}$, and composition is preserved.  However, for
any \emph{particular} monoidal functor $\bV\to \bW$, we could easily
check this directly, without needing to assert the existence of the
whole functor, and thus its very large codomain
$\mathbf{CAT}$.\footnote{One cannot help, however, being reminded of
  how infinite sets in pre-Cantorian mathematics were only regarded as
  `potentialities' rather than completed entities, and how at first
  even ordinary large categories were viewed with suspicion.  It seems
  that the trend of mathematical development is towards recognizing
  ever larger entities as having an independent existence.}  Thus,
often in category theory, the assumption of inaccessibles can be
regarded as merely a convenience (although a very convenient one!).
Thus, it is natural to wonder under what conditions the use of
inaccessibles can be eliminated from categorical arguments.

There is another unsatisfactory aspect of \zfci.  We have asserted
that we only consider objects smaller than \ka\ to be `ordinary
mathematical objects', but what if at some later date we discover that
there happens to be a group larger than \ka\ which we want to include
as a mathematical object?  Clearly this means that we chose the wrong
\ka\ to define `small' and `large' and we should choose a larger one.
In particular, we may want to use category theory to study the
category of categories, and then apply our results to \emph{large}
categories as well as small ones.

To ensure that such switches would always be possible, Grothendieck
proposed an axiom that there are arbitrarily large inaccessibles, or
equivalently that every set is contained in a Grothendieck universe.
(Actually, this was the first axiom using universes to be proposed for
category theory; only later did Mac Lane observe that one universe was
usually sufficient.)  This is still quite a weak large-cardinal axiom.
Note that this use of multiple inaccessibles is different from our
discussion of very large and extremely large categories in
\S\ref{sec:universes}; here we are changing the size of the universe,
rather than using multiple universes at once.

Now, as long as we haven't made use of any properties of \ka\ beyond
its inaccessibility, all our results proven for \ka\ will also be true
for our new, larger, \ka.  This means that all our theorems implicitly
begin with ``For any inaccessible \ka, \dots''.  However, the
arbitrariness of \ka\ may make us somewhat uneasy.  Furthermore, there
is no \emph{a priori} guarantee that the properties of particular
objects will be preserved by change of universe.  For example, suppose
that we prove in \zfci, for some property \ph, that there exists a
small group $G$ such that $\ph(G,H)$ is true for all small groups $H$.
(The assertion that $G$ is a limit or colimit of some specified
diagram in $\mathbf{Grp}$ is of this sort.)  This will then be true
whether `small' is interpreted relative to one inaccessible or
another, but there is no \emph{a priori} reason why the group $G$ with
this property need be the same in the two cases.  Thus we have no way
to conclude that there is a $G$ satisfying $\ph(G,H)$ for \emph{all}
groups $H$.  In particular cases this is obvious; for example, we have
explicit ways to compute (small) limits and colimits in $\mathbf{Grp}$
which do not depend on the size of the universe.  But the absence of a
general truth of this sort means that care is needed when engaging in
such `universe-juggling'.

\section{Natural models and reflection principles}
\label{sec:nmod-refl}

In an attempt to remedy these problems, let us investigate more
specifically what properties of inaccessibility are really necessary
for category theory.  The most important consequence of
inaccessibility of \ka\ appears to be that $V_\ka$ is a model of \zfc\
which is itself a set in our assumed larger model $V$.  Thus it is
natural to look more generally at sets which are models of \zfc.

Now, \emph{a priori} a model of \zfc\ consists only of a set $M$ and a
relation $E\subset M\times M$, to be interpreted as `membership', such
that the axioms of \zfc\ hold.  However, we want the elements of a set
in $M$ to be the same as its elements in $V$, so it is natural to
require that $E$ coincides with the actual membership relation $\in$
in $V$, and that $M$ is \textbf{transitive}, meaning that $x\in y\in
M$ implies $x\in M$.  In fact, any model is isomorphic to a transitive
one, called its \textbf{(Mostowski) transitive collapse}, via an
isomorphism defined inductively by $T(x) = \{T(y) \mid yEx\}$.  Thus,
nothing essential is lost by considering only transitive models.

However, pathologies still exist among transitive models.  In
particular, if there exists any model of \zfc, then there exists one
which is transitive and \emph{countable}---despite the fact that in
\zfc\ one can prove the existence of uncountable sets!  This is known
as \emph{Skolem's paradox}.  The reason is that a set $x$ which is
countable in $V$ may be uncountable to the eyes of $M$, since the
bijection from $x$ to \om\ may not be in $M$.

Skolem's paradox follows from a model-theoretic result called the
\emph{L\"owenheim-Skolem theorem}.  Like G\"odel's incompleteness
theorems and Tarski's undefinability theorem, the L\"owenheim-Skolem
theorem also has important philosophical implications for any
axiomatic foundation of mathematics.  For this reason, and because we
will re-use the same ideas later, I will now sketch a proof of it.

If \ph\ is any statement and $M$ is any structure (that is, a
potential model of some theory, like a group or a universe), let
$\ph^M$ denote \ph\ \textbf{relativized} to $M$, meaning that all its
quantifiers are restricted to range only over elements of $M$.  If
$M\subset N$ we say that $\ph(x_1,\dots,x_n)$ is \textbf{reflected
  from $N$ to $M$} if
\[\forall x_1\in M\cdots \forall x_n\in M \;\Big( \ph^N(x_1,\dots,x_n)
\iff \ph^M(x_1,\dots,x_n)\Big).
\]
If \emph{all} statements \ph\ are reflected from $N$ to $M$, we say
that $M$ is an \textbf{elementary substructure} of $N$ and write $M\ee
N$; clearly in this case $M$ is a model of some theory if and only if
$N$ is.

The \emph{downward L\"owenheim-Skolem theorem} says that given any
structure $N$ and any infinite $S\subset N$, there exists an $M\ee N$
such that $S\subset M$ and $|M| = |S|$.  To construct such an $M$,
start with $M_0=S$.  Now, for each statement $\ph(x_1,\dots,x_n)$ of
the form $\exists y: \psi(y,x_1,\dots,x_n)$ and each $a_1,\dots,a_n
\in M_0$ such that $\ph^N(a_1,\dots,a_n)$ is true, choose some $b\in
N$ such that $\psi^N(b,a_1,\dots,a_n)$ is true.  Let $M_1$ be $M_0$
together with all such `witnesses' $b$.  Iterate this process to
define $M_2,M_3,\dots$, and let $M=\bigcup_{n\in\om} M_n$.  Since
there are only countably many statements \ph, the cardinality never
increases, so $|M| = |S|$.  The only potential difficulty in showing
that $M\ee N$ is with quantifiers; but we have dealt with existential
quantifiers by construction, while universal ones can be rephrased as
the nonexistence of a counterexample (this is called the
\emph{Tarski-Vaught test} for elementary submodels).

We now obtain Skolem's paradox by starting with a model $N$ of \zfc,
applying this theorem for any countably infinite $S\subset N$ (which
exists by the axiom of infinity), then taking the transitive collapse
of the resulting model $M$.  Note that in fact we have proven more: if
\zfc\ (or any theory) has a model $N$, then it has a transitive model
with cardinality \ka\ for all $\ka\le|N|$.  The \emph{upward
  L\"owenheim-Skolem theorem} says that this is also true for $\ka
>|N|$.

To avoid Skolem's paradox, it suffices to require that our transitive
model $M$ of \zfc\ be closed under subsets: $x\subset y\in M$ implies
$x\in M$.  This ensures that if $A,B\in M$ and $A$ and $B$ are
bijective in $V$, then the bijection is also in $M$, since it is a
subset of $A\times B$.  A transitive model of \zfc\ closed
under subsets is called a \textbf{natural model}; see~\cite{mv:nmod}.
It is not difficult to show that any natural model is of the form
$V_\al$ for some limit ordinal \al.

We have seen that $V_\al$ is a natural model when \al\ is
inaccessible, but in fact inaccessibility is much stronger than
necessary.  Inaccessibility of \al\ asserts, in particular, that
$V_\al$ contains the image of any function $f\maps X\to V_\al$ such
that $X\in V_\al$.  But saying that $V_\al$ satisfies the replacement
schema only asserts this when $f$ is \emph{definable} from $V_\al$;
that is, when $f\in\Def(V_\al)$.  Note that in general any such $f$ is
in $V_{\al+1}=\cP V_\al$, which contains $\Def(V_\al)$ as a proper
subset.  Another way of expressing this is to say that $V_\al$ models
\zfc\ if and only if $(V_\al,\Def(V_\al))$ models \nbg, while
$(V_\al,V_{\al+1})$ can only model \nbg\ when \al\ is inaccessible (in
which case it also models \mk).  To see how nontrivial this
distinction is, observe that since there are only countably many
statements \ph, we have $|\Def(V_\al)| = |V_\al|$, while of course
$|V_{\al+1}| = 2^{|V_\al|} > |V_\al|$.


\begin{rmk}
  One can also say that if \al\ is inaccessible, then $V_\al$
  satisfies the \emph{second-order replacement axiom}.  To understand
  this we need to describe second-order logic.  The logic we have
  discussed so far is called \emph{first-order logic} because
  variables and quantifiers only range over `things'; in second-order
  logic they are also allowed to range over `sets of things'.  This
  distinction can be confusing, since for a first-order theory like
  \zfc\ the `things' are themselves called `sets'.

  Second-order logic is more powerful than first-order logic, but
  suffers from an ambiguity of interpretation.  If by `set of things'
  we intend to mean \emph{any} subset of the model under
  consideration, then we need an external set theory to define what is
  meant by this.  This is the sense in which a Grothendieck universe
  satisfies the `second-order replacement axiom'.  On the other hand,
  if we allow the model itself to stipulate what `sets of things'
  exist, then second-order logic reduces to first-order logic
  augmented with an extra sort of `thing' called a `set of things'.
  This occurs, for example, with the replacement axiom of \nbg.  The
  existence of natural models that are not Grothendieck universes
  shows that this ambiguity has teeth.  For this reason I will
  continue to stick to first-order logic.
\end{rmk}

Now, by slightly modifying the proof of the L\"owenheim-Skolem
theorem, we can prove that for any \be\ and $S\in V_\be$, then there
is an $\al\le \be$ such that $S\in V_\al$ and $V_\al\ee V_\be$.
Namely, instead of a sequence $M_0\subset M_1\subset M_2\subset\dots$,
we construct a sequence $\al_0<\al_1<\al_2<\dots$, by letting
$\al_{n+1}$ be some ordinal such that all the witnesses $b$ for
$V_{\al_n}$ are contained in $V_{\al_{n+1}}$.  Setting $\al =
\lim_{n<\om} \al_n$, it follows that $V_\al\ee V_\be$.  In particular,
if $V_\be$ is a natural model, so is $V_\al$.

Of course, in this case there is no guarantee that $\al\neq\be$.  But
if \be\ is \emph{inaccessible}, then at the $n^{\mathrm{th}}$ stage we
need to add only $|V_{\al_n}|<\be$ witnesses, so we can choose
$\al_{n+1}<\be$; and since \be\ is regular, we then have $\al<\be$.  A
slight improvement of this proof shows that if \be\ is inaccessible,
then
\[\{\al<\be \mid V_\al \text{ is a natural model}\}
\]
is stationary in \be, and in particular has cardinality \be.  Thus,
the existence of natural models can be regarded as a sort of `large
cardinal axiom' significantly weaker than even a single inaccessible
(although \al\ need not actually be a cardinal for $V_\al$ to be a
natural model).

It makes intuitive sense that we could also carry through the above
argument using the whole universe $V$ instead of $V_\be$, and thereby
construct a natural model $V_\al$ with $V_\al\ee V$.  Of course, this
would prove Con(\zfc) and violate the incompleteness theorem.  The
flaw in the argument is that with $V$ in place of $V_\be$, Tarski's
undefinability theorem prevents us from defining the set of witnesses
$b$.  However, we can rescue the argument if instead of asserting
$V_\al\ee V$, we only assert that a given \emph{finite} set of
statements is reflected from $V$ to $V_\al$.  This gives a
theorem-schema called the \emph{reflection principle}: for any finite
set $\ph_1,\dots,\ph_n$ of statements, we can prove in \zfc\ that for
any set $y$, there exists an \al\ such that $y\in V_\al$ and
$\ph_1,\dots,\ph_n$ are all reflected in $V_\al$.

Now, while the single statement ``there exists a natural model''
implies the consistency of \zfc, the reflection principle suggests a
version of it that does not.  Let \bbS\ be an extra constant symbol,
and add to \zfc\ the axiom ``\bbS\ is transitive and closed under
subsets'' along with a \emph{reflection axiom}
\[\forall x_1\in \bbS\cdots \forall x_n\in \bbS \;\Big( \ph(x_1,\dots,x_n)
\iff \ph^\bbS(x_1,\dots,x_n)\Big).
\]
for each statement \ph.  We denote the resulting system by by
\textbf{ZFC/S} (``\zfc\ with smallness'').  It follows that for each
axiom \ph\ of \zfc, the relativized version $\ph^\bbS$ is true in
\zfcs, so \bbS\ is a model of \zfc.  We can also show $\bbS=V_\bbk$
for some ordinal \bbk.  However, since the proof of each axiom
$\ph^\bbS$ uses a different instance of the reflection axiom schema,
we cannot prove in \zfcs\ the single statement ``\bbS\ is a natural
model''.

In fact, we can prove that \zfcs\ is \emph{conservative} over \zfc.
Suppose we have any theorem which is provable in \zfcs; we show that
it can also be proven in \zfc.  The original proof, being only
finitely long, can only use the reflection schema for finitely many
statements \ph.  Thus, by the reflection principle for \zfc, we can
find an \al\ such that all these statements \ph\ are reflected in
$V_\al$.  We can then replace $\bbS$ by $V_\al$ and carry out the
proof in \zfc.  Note, however, that unlike the proof for
conservativity of \nbg, we do not have a way to construct, from a
model of \zfc, a model of \zfcs\ with the same (small) sets.

The theory \zfcs\ is due to Feferman, who proposed it
in~\cite{feferman:fdns-of-ct} as a foundation for category theory.  Of
course, we now define \textbf{small} to mean ``element of \bbS'' and
\textbf{large} to mean ``set not necessarily in \bbS''; otherwise
things go mostly as before.  Since the axioms of \zfc\ are satisfied
for all sets, we can manipulate large sets as we wish; so we retain
that advantage of \zfci.  If we want to talk about very large and
extremely large objects, it is easy to add multiple symbols
$\bbS_1,\bbS_2,\dots$, each satisfying reflection and with
$\bbS_1\in\bbS_2\in\cdots$.

However, because \zfcs\ is conservative over \zfc, we have not
strengthened our basic set theory.  In particular, anything about
small objects that we prove with the aid of large categories would
still be provable in pure \zfc.  Thus, we obtain a precise version of
our intuition that the use of inaccessibles in category theory is
merely for convenience: since many categorical proofs stated using
inaccessibles can be formalized in \zfcs, any \emph{consequence} of
such a theorem not referring explicitly to inaccessibles is also
provable purely in \zfc.

\zfcs\ also eliminates at least some sources of universe-juggling.
For example, because every statement about sets is reflected in \bbS,
anything we prove in \zfcs\ about small objects is \emph{also} true
about \emph{large} objects.  In particular, anything we prove about
small categories, even making use of the large category
$\mathbf{Cat}$, will also be true about large categories.  Moreover,
any property of small objects which refers only to small objects is
retained when reinterpreted to refer to large objects.  For example,
the statement that $G$ is a limit of a specified small diagram in the
category $\mathbf{Grp}$ of small groups can be expressed as
$\ph^\bbS(G)$ for some statement \ph.  Thus, the reflection principle
implies that $\ph(G)$ is also true, and hence $G$ satisfies the same
universal property with respect to large groups.

Whether \emph{all} universe-juggling can be eliminated in this way
depends on our ontological position towards \zfcs.  If we believe that
\zfcs\ is a true representation of reality---that is, that there
actually \emph{exists} an \bbS\ satisfying the reflection axioms, and
when working in \zfcs\ we are making statements about that particular
\bbS---then of course not all objects are small.  The reflection
principle gives us a precise sense in which they `might as well be'
small, but if we insist on being able to make them \emph{actually}
small, we would need to augment \zfcs\ by a Grothendieck-like
assumption of many natural models and continue to engage in
universe-juggling.

However, if we instead take the position that \zfc\ is `true', while
\zfcs\ is only a convenient fiction made possible by the reflection
principle, then we can relegate all the universe-juggling to the
`behind the scenes' interpretation of \zfcs\ in \zfc.  That is, we
prove all our results in \zfcs, and argue that when we want to apply
them to particular objects in the `real world' of \zfc, we tacitly use
the reflection principle to choose some $V_\al$ which contains all the
objects we happen to be interested in.  Thus the theorems of category
theory take on the character of a `meta-theory', which can be applied
to any particular set of objects in the real world by choosing a
sufficiently large $V_\al$ containing that set.

\begin{rmk}
  This seems an appropriate place to mention \textbf{Ackermann set
    theory}, a theory of sets and classes like \nbg\ and \mk\ which
  has also been proposed as a foundation for category theory.  Unlike
  \nbg\ and \mk, however, it allows classes to be elements of other
  classes.\footnote{This is not obvious from the axioms, which only
    assert directly the existence of classes whose elements are sets.
    To see that there must be classes containing other classes, we
    observe first that axiom~\ref{item:ack4} implies that sethood
    cannot be characterized without referring to it explicitly;
    otherwise the class $V$ of all sets would be a set.  This means
    that since the property ``$\exists y: x\in y$'' is true of all
    sets, it must also be true of some classes; otherwise it would
    characterize sethood.} Its axioms are the following.
  \begin{enumerate}
  \item Extensionality, foundation, and choice.
  \item Any element or subset of a set is a set.
  \item For any definable property $\ph(x)$, if all objects $x$
    satisfying $\ph(x)$ are sets, then the class $\{x\mid \ph(x)\}$
    exists.
  \item If in the previous axiom, in addition $\ph(x)$ does not refer
    explicitly to sethood (that is, to whether or not any given class is
    a set), then the class $\{x\mid \ph(x)\}$ is a set.\label{item:ack4}
  \end{enumerate}
  These axioms imply (though not obviously) that the class $V$ of all
  sets is a model of \zfc.  Conversely, any model of \zfcs\ satisfies
  these axioms if `class' means `set' and `set' means `element of
  \bbS'; see~\cite{levy:ackermann,reinhardt:ackermann}.  Thus, like
  \nbg\ and \zfcs, Ackermann set theory is a conservative extension of
  \zfc.\footnote{The intuition behind Ackermann set theory, however,
    is different from that of \zfc.  Ackermann argued that that the
    elements of a \emph{set} must be `sharply delimited', while the
    elements of a \emph{class}, such as $V$, may depend on how broadly
    we interpret the concept of `set'.  Thus, only properties which do
    not refer explicitly to sethood are `sharply delimited' enough to
    define sets.  It is striking that nevertheless, Ackermann's axioms
    turned out \emph{a posteriori} to be equivalent to \zfc\ in what
    they can say about sets.}  But while it implies a limited
  reflection principle, overall it is strictly weaker than \zfcs\ in
  what it can say about its classes.  Most of the subsequent remarks
  about \zfcs\ apply just as well to Ackermann set theory.
\end{rmk}

On the other hand, \zfcs\ (and likewise Ackermann set theory) is not
quite the paradise it first appears.  In order to enable ourselves to
manipulate large objects freely without strengthening set theory, we
have been forced to weaken the replacement axiom for small sets.  In
both \nbg\ and \zfci, we have a replacement axiom saying, essentially,
that the image of a small set under a large function is small.  In
\zfcs, however, we can only assert this if the large function is
\emph{small-definable} (that is, in $\Def(\bbS)$).  This distinction
is invisible to \zfc\ and \nbg\ because there, all functions are (or
might as well be) small-definable, while it is invisible to \mk\ and
\zfci\ because their stronger axioms guarantee that all functions with
small domain are small, even those that are not \emph{a priori}
small-definable.

Perhaps surprisingly, it turns out that this weakening of replacement
has significantly annoying consequences for category theory.  For
example, it implies that the category $\Set=\Set[\bbS]$ of small sets
need not admit limits and colimits for all functors $F\maps \bA\to
\Set$ when \bA\ is small, but only those for which \emph{$F$ is also
  small}.  The same is true for other large categories constructed
from \Set.  Similarly, for a functor $u\maps \bA\to\bB$ between small
categories, the functor
\[\pmb[\bB,\Set\pmb] \xto{-\circ u} \pmb[\bA,\Set\pmb]
\]
will not in general have left and right adjoints (Kan extensions).

On the surface, these restrictions appear quite problematic; the
completeness and cocompleteness of \Set\ is certainly of central
importance in category theory.  One way to respond is to assert that
in \zfcs, the correct definition of \emph{complete} is ``having limits
for all small functors'' rather than ``having limits for all functors
with small domain''.  Similarly, if we write $\lbb\bA,\Set\rbb$ for
the full subcategory of $\pmb[\bA,\Set\pmb]$ determined by the
\emph{small} functors, then induced functor
\[\lbb\bB,\Set\rbb \xto{-\circ u} \lbb\bA,\Set\rbb
\]
does have both adjoints, and we may assert that in \zfcs, the correct
presheaf category to consider is $\lbb\bA,\Set\rbb$ rather than
$\pmb[\bA,\Set\pmb]$.

Another advantage of $\lbb\bA,\Set\rbb$ is that unlike
$\pmb[\bA,\Set\pmb]$, it is small-definable.  This is a good thing,
because in \zfcs, only small-definable categories are at all
well-behaved.  For instance, let $A$ be a small set of objects in a
locally small and small-definable category \bB.  We can then prove
that
\begin{enumerate}[(a)]
\item the full subcategory \bA\ of \bB\ determined by the objects in
  $A$ is small;
\item the inclusion functor $i\maps \bA\to\bB$ is small;
\item for any object $X\in \bB$, the restricted hom-functor
  $\bA(i-,X)\maps \bA\op\to\Set$ is small; and thus
\item the restricted Yoneda embedding $\bB\to \pmb[\bA\op,\Set\pmb]$
  factors through $\lbb\bA\op,\Set\rbb$.
\end{enumerate}
However, there seems to be no way to prove any of these statements if
\bB\ is not small-definable.  Similarly, theorems like the Adjoint
Functor Theorem only seem to work for small-definable categories.

As observed by Feferman, all this makes little difference in most
concrete applications, because any particular diagram, functor, or
category we are interested in will generally be small-definable (at
least, up to equivalence).  However, this is not always trivial to
verify; we saw in \S\ref{sec:universes} examples of categories that
were equivalent to small-definable ones, but not obviously.  Moreover,
small-definability restrictions are tiresome to keep track of, and
some would say unaesthetic as well.  In the next two sections, we will
explore two ways to deal with this problem.



\section{Strong reflection principles}
\label{sec:strong-reflect}

There is an obvious way to `have our cake and eat it too': we can add
to \zfcs\ the extra axiom ``\bbS\ is a Grothendieck universe'', or
equivalently $\bbS=V_\bbk$ where \bbk\ is inaccessible.  For reasons
to be explained below, I will call the resulting theory
\textbf{ZMC/S}.  It then follows, as in \zfci, that every functor with
small domain is small, and all small-definability restrictions vanish.

Of course, this may seem like a step backwards, since we began the
previous section by looking for a way to \emph{avoid} inaccessibles.
However, along the way we discovered the reflection principle, and we
saw that the reflection axiom-schema of \zfcs\ is really what resolves
many of the problems with \zfci\ and allows us to avoid
universe-juggling.  Since \zmcs\ retains reflection, all of this is
still true, so the only disadvantage of \zfci\ which carries over to
\zmcs\ is that it is not conservative over \zfc.

In fact, \zmcs\ is significantly stronger than \zfci, since reflection
implies that \bbk\ is far from the smallest inaccessible.  Namely,
since there exists an inaccessible, there must exist a small
inaccessible; but then there exist two inaccessibles, and so there
must exist two small inaccessibles, and so on.  By applying reflection
to the statement ``there exists an inaccessible larger than \al'', in
\zmcs\ we can even derive Grothendieck's axiom that there are
arbitrarily large inaccessibles.

The same argument that shows \zfcs\ to be conservative over \zfc\
shows that \zmcs\ is conservative over \zfc\ + ``any finite set of
formulas is reflected in some Grothendieck universe.''  This stronger
reflection principle turns out to be equivalent to the assertion that
any closed unbounded (definable) class of ordinals contains an
inaccessible, which essentially says that the cardinal of the universe
is Mahlo (see~\cite{levy:strong-inf}).  \zfc\ augmented by this axiom
is sometimes called \zmc\ (whence the notation \zmcs).  While rather
stronger than the existence of a single inaccessible, and stronger
even than Grothendieck's axiom, this is still quite weak compared to
many large-cardinal axioms, as we saw in \S\ref{sec:large-cardinals}.

Moreover, the principle incorporated in \zmc\ is one of the most
easily motivated large-cardinal axioms.  For instance, it can be
argued that it is a straightforward expression of the
`inexhaustibility' of the universe of sets by any finite number of
operations.  Additionally, it is not difficult to show
(see~\cite{levy:strong-inf}) that in the presence of the basic axioms
only, the reflection principle of \zfc\ implies replacement, power
set, and infinity---all the axioms of \zfc\ which produce larger and
larger sets.  It follows that \zmc\ is equivalent to just the basic
axioms and choice together with the schema ``any finite set of
formulas is reflected in some Grothendieck universe.''  I find this
aesthetically quite appealing, because it captures exactly what
category theory seems to need from set theory: we may not be able to
have a category of \emph{all} sets, but for any particular purpose, we
can choose a category of sets large enough that it might as well
contain all of them.

\section{Toposes and indexed categories}
\label{sec:topos-indexed}

There is another way to deal with the small-definability issues in
\zfcs: we can use \emph{indexed categories}, a tool developed to solve
a similar problem in elementary topos theory.  Since topos theory is
of interest in its own right as a foundation for mathematics in
general, and category theory in particular, we start with a summary of
it.  Good introductions to topos theory can be found
in~\cite{mm:shv-gl,mclarty:ecat-etop}, while~\cite{ptj:elephant} is
encyclopedic.

By way of motivation, observe that while \zfc\ suffices as a
foundation for most of mathematics, there is a sense in which it is
disconnected from most mathematical practice.  In \zfc\ there is a
\emph{global} membership predicate, meaning that if I give you two
random sets, it makes sense to ask whether one is an element of the
other.  However, in actual mathematical practice we usually only speak
of \emph{local} set membership, meaning that asking whether $x\in A$
is only meaningful in the context of some fixed set $B$ such that $x$
is known to be an element of $B$ and $A$ is known to be a subset of
$B$.  In other words, the way most mathematicians usually think of a
set (or a group, or a topological space, etc.) is as a collection of
`abstract' elements which have no `internal' structure aside from
being elements of that set.  The only way that elements of two
different sets relate to each other is via functions and relations
between those sets.

Of course, this is a very categorical way of thinking; it is closely
related to the assertion that we only care about objects of a
category, such as \Set, up to isomorphism.  Thus, it is natural to try
to axiomatize the properties of the category \Set, instead of
axiomatizing a global membership relation.  (While useful, this
motivation for topos theory is ahistorical;
see~\cite{mclarty:hist-topos}.)  The appropriate axioms can be
classified just as we did for \zfc\ and \nbg.
\begin{enumerate}
\item The \emph{basic} axioms: \Set\ is cartesian closed, has finite
  limits and colimits, and the terminal object is a generator.
\item The \emph{size-increasing} axiom: every object has a `power
  object' classifying its subobjects.
\item The \emph{size-assertion} axiom: there is a `natural numbers
  object'.
\item The axiom of \emph{choice}: all epimorphisms split.
\end{enumerate}
An \textbf{elementary topos}\footnote{There is some controversy about
  whether the English plural of \emph{topos} should be \emph{toposes}
  or \emph{topoi}.  The paper by Grothendieck and Verdier which coined
  the term is in French, where the plural is again \emph{topos}; this
  seems to tell against the Greek plural.  However, in English
  \emph{topos} also means ``a literary theme or motif'', and in this
  case the plural used is always \emph{topoi}.  Unfortunately one
  cannot avoid making one choice or the other!}  is a category with
finite limits and power objects; this implies cartesian closedness and
the existence of finite colimits.  A \textbf{natural numbers object
  (NNO)} in a topos is an object satisfying the universal property of
definition by recursion, or equivalently proof by induction.  A topos
is \textbf{well-pointed} if the terminal object $1$ is a generator;
this means that an object $X$ is determined by its `elements' $x\maps
1\to X$.  Thus, a natural axiomatization of \Set\ is that it is a
well-pointed elementary topos with a \nno\ and satisfying the axiom of
choice (a \textbf{WPTNC}).  These axioms for \Set\ are also referred
to as the \textbf{Elementary Theory of the Category of Sets (ETCS)},
after Lawvere's influential
paper~\cite{lawvere:etcs,lawvere:etcs-long}.  Just as with \zfc, it is
an empirical observation that much of mathematics can be developed
starting from a model of \etcs.

\begin{rmk}\label{rmk:constructive-1}
  In fact, much mathematics can be developed from any elementary topos
  (perhaps with a \nno), as long as one uses \emph{intuitionistic} or
  \emph{constructive logic} instead of classical logic.  I will return
  to this point in \autoref{rmk:constructive-2}.
\end{rmk}

Notably absent from \etcs\ is any analogue of the axiom of
replacement.  This means that if \al\ is any limit ordinal greater
than \om, so that $V_\al$ satisfies all the axioms of \zfc\ except
replacement, then the category $\Set[V_\al]$ of sets and functions in
$V_\al$ is a \wptnc.  As remarked in \S\ref{sec:ordinals}, most
mathematical objects outside of set theory have very low rank, living
quite comfortably in $V_{\om\cdot 2}$, so this lends extra credence to
the observation that \etcs\ suffices for much of mathematics.

In fact, all one needs to construct a \wptnc\ is a model of
\textbf{Bounded Zermelo set theory with Choice (BZC)}, which is \zfc\
without replacement and in which the properties in the separation
axiom are only allowed to have bounded quantifiers (``for all $x\in
A$'' rather than just ``for all $x$'').  This version of separation is
variously called \textbf{restricted}, \textbf{bounded}, or
\textbf{$\Delta_0$-separation}.  Conversely, from a \wptnc\ one can
construct a model of \bzc, although some cleverness is needed to
obtain a global membership predicate; see~\cite[VI.10]{mm:shv-gl}
or~\cite[Ch.~9]{ptj:topos-theory} for two approaches.  Thus, \etcs\
and \bzc\ are \emph{equiconsistent}.

I will discuss the implications of \etcs\ and \bzc\ for mathematical
practice in more detail in \S\ref{sec:cst}; for now, let us consider
their consequences for category theory.  We have seen a hint already
of what can go wrong without replacement in our study of \zfcs, where
weakening the replacement axiom created unexpected problems.  To see
the more drastic problem we are now faced with, consider the meaning
of the statement ``\bA\ has small products'' when \bA\ is a large
category.  Intuitively, this means that any $X$-indexed family of
objects of \bA\ has a product, for any `set' $X$.  But what exactly is
an ``$X$-indexed family of objects of \bA''?

In \nbg\ this can mean a class which is a function from $X$ to the
class of objects of \bA, while in \zfc\ it can mean a definable
property \ph\ such that for any $x\in X$ there is a unique $a\in \bA$
with $\ph(x,a)$.  In either case, we can then apply the replacement
axiom to obtain a \emph{set} such as $\{a\in \bA \mid \exists x\in X\;
f(x)=a\}$, and go on to construct the product.  However, if $X$ is
instead an object of an elementary topos \bS, it is not at all clear
what is even \emph{meant} by an $X$-indexed family of objects of \bA.

On the other hand, any large category \bA\ which is \emph{constructed
  from} \bS\ in any reasonable way will come with a canonical notion
of $X$-indexed family of objects for any $X\in\bS$.  To start with,
\bS\ itself has such a notion: an $X$-indexed family of objects of
\bS\ can be defined to be simply an arrow $p\maps K\to X$.  The
intuition is that an element $x\in X$ indexes the fiber $p^{-1}(x)$,
but of course objects of the abstract category \bS\ have no `elements'
as such.  It follows that the \emph{category of} $X$-indexed families
of objects of \bS\ should be the slice category $\bS/X$.  We can then
extend this to other categories of `sets with structure'.  For
example, if \bA\ is the category of `small groups', meaning internal
group objects in \bS, then an $X$-indexed family of objects of \bA\
should be an internal group object in $\bS/X$.

Thus, when doing category theory relative to an elementary topos \bS,
it is natural to consider, instead of a single category \bA, a family
of categories $\bA^X$ for each object $X\in\bS$, where $\bA^X$ is
thought of as the category of `$X$-indexed families of objects of
\bA'.  Of course, these categories should be related in some way as
$X$ varies, and it turns out that the important property is the
ability to \emph{reindex} a family: if $K$ is an $X$-indexed family
and $f\maps Y\to X$, then we have a $Y$-indexed family $f^*K$.  The
intuition is that for $y\in Y$ we have $(f^*K)_y=K_{f(y)}$.  For \bS\
and categories constructed from it, this reindexing is given by
pullback along $f$.

With this as motivation, we define an \textbf{\bS-indexed category}
\bA\ to consist of a category $\bA^X$ for each object $X\in\bS$,
together with a functor $f^*\maps \bA^X\to\bA^Y$ for each arrow
$f\maps Y\to X$ in \bS, and natural isomorphisms $(gf)^*\iso f^*g^*$
and $\Id_{\bA^X} \iso (1_X)^*$ satisfying obvious axioms.\footnote{The
  2-categorically sophisticated reader may call this a pseudofunctor
  $\bS\op\to\mathcal{CAT}$.  This is fine as long as we have some
  external set theory with which to define a 2-category
  $\mathcal{CAT}$ of large enough categories.}  An \bS-indexed
category can also be described by assembling the categories $\bA^X$
into a single category $\oint\!\bA$ equipped with a functor
$\oint\!\bA\to \bS$ assigning a family to its indexing object; in this
form it is called a (categorical) \textbf{fibration} over \bS.  Good
introductions to indexed categories and fibrations can be found
in~\cite[Part~B]{ptj:elephant} and~\cite{streicher:fibredcats}.

If we now replace our na\"ive large category \bA\ by an \bS-indexed
category, all problems disappear.  For example, we can define \bA\ to
have \emph{\bS-indexed products} if each reindexing functor $f^*\maps
\bA^X\to\bA^Y$ has a right adjoint (plus a commutativity condition).
For a more general notion of completeness, we need a notion of `small
category', and the obvious candidate is an \emph{internal category} in
\bS, which consists of objects $C_0,C_1\in\bS$ with arrows $s,t\maps
C_1\toto C_0$, $i\maps C_0\to C_1$, and $c\maps C_1\times_{C_0} C_1\to
C_1$ satisfying obvious axioms.  Any internal category $C$ gives rise
to an \bS-indexed category \bC\ with $\bC^X = \bS(X,C_\bullet)$, and
we can define a \emph{$C$-diagram} in any \bS-indexed category \bA\ to
be an object $F\in \bA^{C_0}$ together with a morphism $s^* F\to t^*F$
in $\bA^{C_1}$ satisfying suitable axioms.  The appropriate notions of
\emph{limit} and \emph{completeness} are then fairly
straightforward.\footnote{This is not quite true; I am omitting some
  details in an attempt to give the flavor of the subject without
  getting bogged down.  See the references above for a careful
  treatment.}  Similarly, we can define local smallness, generators,
and well-poweredness, and state and prove an Indexed Adjoint Functor
Theorem.

We saw that \bS\ itself is represented by the \bS-indexed category
with $\bS^X = \bS/X$; we call this the \textbf{self-indexing} of \bS.
It turns out that the self-indexing is always complete and cocomplete,
in the indexed sense, for \emph{any} elementary topos \bS.  In
particular, this applies to $\Set[V_\al]$ for any limit ordinal \al,
even though such toposes are not complete and cocomplete in the
na\"ive `external' sense.  For example, $\Set[V_{\om\cdot 2}]$ fails
to have the coproduct $\coprod_{n<\om}\beth_{\om+n}$, even though it
contains $\om$ and each $\beth_{\om+n}$.  But this does not violate
indexed cocompleteness, since $\big\{\beth_{\om+n}\big\}_{n\in\om}$ is
not an $\om$-indexed family in the self-indexing of $\Set[V_{\om\cdot
  2}]$: if there were such a family $K\to \om$, then $K$ would have to
essentially already \emph{be} the desired coproduct.

The situation may be clarified by observing that there is a
\emph{different} indexed category over any $\bS=\Set[V_\al]$ called
the \textbf{na\"ive indexing}, for which $\bS_{\text{na\"ive}}^X$ is
the category of all functions $X\to V_\al$.  The self-indexing of
$\Set[V_\al]$ embeds in its na\"ive indexing, and the two are
equivalent precisely when \al\ is inaccessible.  In between the two we
have the \textbf{definable indexing}, for which $\bS_{\text{def}}^X$
is the category of all \emph{definable} functions $X\to V_\al$.  This
agrees with the self-indexing whenever $V_\al$ is a natural model.
Neither the na\"ive indexing nor the definable indexing of
$\Set[V_{\om\cdot 2}]$ is complete or cocomplete, nor is the na\"ive
indexing of $\Set[V_\al]$ when $V_\al$ is a natural model that is not
a Grothendieck universe; but the self-indexing of either always is.

In fact, if $V$ is \emph{any} model of \bzc, then the resulting
well-pointed topos $\Set[V]$ has both a self-indexing and a definable
indexing, and the assertion that $V$ satisfies replacement (hence is a
model of \zfc) is equivalent to the assertion that these two indexings
agree.  To define a `na\"ive indexing', our universe $V$ must live
inside a larger universe of sets or classes, such as a model of \nbg,
\mk, or \zfci.  In each of these cases, the version of replacement
asserted implies that the resulting na\"ive indexing is actually
equivalent to the self-indexing as well.

We can now see that when working with large categories in \zfc, we
have implicitly been using the definable indexing, while in \nbg, \mk,
and \zfci\ we have been using the na\"ive indexing.  Topos theory
suggests that actually, the \emph{self-indexing} is always the
`correct' indexing to use, and the role of the replacement axiom is to
ensure that the self-indexing is equivalent to the definable or
na\"ive indexings, which are more intuitive and easier to work with.
See~\cite[\S 17]{streicher:fibredcats} for further discussion of
`internal' versus `external' completeness.

Now let us return to \zfcs.  Since \bbS\ is a natural
model\footnote{However, remember that ``\bbS\ is a natural model'' is
  not a single theorem of \zfcs, but a schema consisting of one
  theorem for each axiom of \zfc.} but not a Grothendieck universe,
the self-indexing of $\Set[\bbS]$ agrees with its definable indexing,
but not its na\"ive indexing.  Moreover, the objects of
$\bS^X_{\text{self}} = \bS^X_{\text{def}}=\bS/X$ are essentially the
same as the \emph{small} functions $X\to\bbS$; thus indexed
completeness of $\Set[\bbS]$ agrees with our proposed \emph{ad hoc}
redefinition of completeness in \S\ref{sec:nmod-refl}.  More
generally, any small-definable category in \zfcs\ gives rise to a
$\Set[\bbS]$-indexed category containing only the \emph{small}
$X$-indexed families of objects, and the machinery of indexed
categories automatically keeps track of all the restrictions we had to
impose by hand in \S\ref{sec:nmod-refl}.  For example, if $C$ is an
internal category $C$ in $\Set[\bbS]$, the category of $C$-diagrams in
$\bS_{\text{self}}$ is the well-behaved $\lbb C,\Set\rbb$ rather
than the poorly-behaved $\pmb[C,\Set\pmb]$.

Thus, one may say that the problems we encountered with \zfcs\ arose
due to our trying to use the na\"ive indexing in a situation where our
replacement axiom was only sufficient to deal with the definable
indexing.  Moreover, Feferman's hypothesis that \zfcs\ suffices for
basic category theory now follows from the observation that most
theorems of basic category theory have indexed analogues.

On the other hand, once we are willing to use indexed categories, we
do not necessarily need to to assume a replacement axiom at all in
order to do category theory.  Thus, we could choose \emph{any} \wptnc\
\bS\ which is a set, even $\Set[V_{\om\cdot 2}]$, define a \emph{small
  category} to mean an internal category in \bS\ and a \emph{large
  category} to mean an \bS-indexed category, and develop category
theory that way.  This way we can remain completely within \zfc\ (or
even something much weaker), and Grothendieck's axiom of arbitrarily
large inaccessibles can be replaced by the simple fact that any set is
contained in some $V_\al$.

However, the machinery of indexed categories is admittedly rather
complicated, and it seems unreasonable to expect most users of
category theory to be familiar with it when there are so many simpler
foundational options available.  If nothing else, though, indexed
categories give a more conceptual understanding of the
small-definability restrictions arising in \zfcs.  (Of course, indexed
categories are also crucial when working with a general elementary
topos, rather than a \wptnc.)

\section{Aside: the strength of categorical set theory}
\label{sec:cst}

Let us pause here briefly to compare the `category-theoretic
foundation' for mathematics offered by \etcs\ and its relatives with
the `set-theoretic foundation' offered by \zfc\ and its cousins.  This
terminology is common, but one can also argue persuasively
(see~\cite{lawvere:etcs-long}) that \etcs\ is itself a \emph{set
  theory}, meaning a theory about the behavior of sets.  What
distinguishes it from \zfc\ is not its objects of study, but how it
studies them: by taking functions as a basic notion rather than global
membership.  Perhaps a more correct distinction would be to call
\etcs\ a \emph{categorical\footnote{There is an unfortunate collision
    between the common and natural use of `categorical' to mean
    `related to categories', and the much older philosophical and
    logical tradition in which `categorical' means `absolute' or
    `uniquely determined'.  This has led some authors to use
    \emph{categorial} for the former notion.} set theory} and \zfc\ a
\emph{membership set theory}.

I mentioned in \S\ref{sec:topos-indexed} that \etcs\ is equiconsistent
with \bzc.  In fact, if we add axioms to \etcs\ and \bzc\ saying that
every set is contained in a transitive one and that transitive
collapses exist (\S\ref{sec:nmod-refl}), then we can obtain a full
equivalence between models of the two theories;
see~\cite[Ch.~9]{ptj:topos-theory} or~\cite{osius:cat-setth}.  These
additional axioms can be proven in \zfc\ using replacement, but are
much weaker than it; in particular, adding them does not change the
consistency strength of \etcs\ and \bzc.  This altered version of
\bzc\ is sometimes called \textbf{Mac Lane set theory (MAC)}.  Thus,
at least in one sense, \etcs\ and \bzc/\mac\ are completely
equivalent.

However, as I argued in \S\ref{sec:topos-indexed}, \etcs\ may seem
closer to most mathematical practice than \bzc, since it discards the
usually superfluous notion of global membership.  Furthermore, in line
with the observed fact that most mathematicians only care about
objects up to isomorphism, \etcs\ can only \emph{characterize} any set
up to isomorphism; see~\cite{mclarty:numbers}.  Even a lot of notions
in set theory, including many large-cardinal axioms, are invariant
under isomorphism.

On the other hand, \etcs\ and \bzc\ are both significantly weaker than
\zfc: not only are they missing the replacement axiom, but they only
allow separation for formulas with bounded quantification.  This
implies that just as \nbg\ can only prove mathematical induction for
statements not quantifying over classes, \bzc\ can only prove
induction for statements without unbounded quantifiers.  For example,
if \bA\ is a large category in the style of \S\ref{sec:classes}, then
a statement such as ``for all $n$, \bA\ has $n$-fold products''
involves unbounded quantifiers and thus cannot be proven by induction
in \bzc\ or \etcs, at least not obviously.

The axiom of replacement also has other uses outside of set theory,
usually taking the form of transfinite induction arguments.  The
classic example is \emph{Borel determinacy} in descriptive set theory,
which is known to be unprovable even in Zermelo set theory (that is,
\zfc\ without replacement but with full separation).  Closer to home
for us are transfinite constructions such as \autoref{thm:soa}, which
also require some form of replacement.  For instance, the power-set
functor cannot be iterated even \om\ times without replacement, since
$|\cP^\om \om|=\beth_\om \not\in V_{\om\cdot 2}$.  The same is true of
the dual-vector-space functor.  For a very detailed study of the
strength of \bzc\ and its cousins, see~\cite{mathias:str-maclane}.

For these reasons, it is natural to wonder whether \etcs\ can be
strengthened with versions of full separation and replacement to
obtain a categorical set theory equivalent in strength to \zfc.  In
fact, it suffices to consider replacement, since at least with
classical logic, replacement implies full separation.  There has not
been a lot of work in this area, but several sorts of categorical
replacement-type axioms have been proposed.

The perspective of \S\ref{sec:topos-indexed} suggests that perhaps a
categorical replacement axiom should essentially say ``the definable
indexing of \Set\ agrees with the self-indexing''.  In constructing
the definable indexing of $\Set[V]$ we used the fact that our sets
have elements, but in a well-pointed topos we can replace these
elements by morphisms $1\to X$.  Thus, we say a well-pointed topos
\bS\ \emph{satisfies replacement} if for any object $X$ and any
definable property $\ph(x,S)$ such that for any `element' $x\maps 1\to
X$ there exists a object $S_x$ unique up to isomorphism with
$\ph(x,S_x)$, there exists a morphism $S\to X$ such that for any $x$
there is a pullback square
\[\xymatrix{S_x \ar[r]\ar[d] & S \ar[d]\\
  1 \ar[r]_x & X.}\]
This version of replacement is from~\cite{mclarty:catstruct}; other
similar ones can be found in~\cite{osius:cat-setth,lawvere:etcs-long}.
These references also prove that \etcs\ plus replacement is equivalent
to \zfc, in the strong sense of an equivalence of models.  Thus, if
all we want is a categorical set theory equivalent to \zfc, we have
it.


On the other hand, these axioms of replacement are not fully
satisfactory from a categorical point of view, because they all depend
heavily on well-pointedness.  As mentioned previously, all the other
axioms of \etcs\ make perfect sense without well-pointedness, and much
mathematics can be developed from any elementary topos; thus it would
be nice to have a version of the replacement axiom that makes sense
without well-pointedness.

In~\cite[Ch.~9]{taylor:pracfdn} Paul Taylor proposed an axiom he
called the \emph{categorical axiom of iterative replacement}, which
asserts directly the possibility of transfinite constructions on
functors.  This axiom makes sense without well-pointedness, but
seemingly has not been investigated in very much detail.  I do not
know of any more explicitly replacement-like axioms that make sense
for non-well-pointed toposes.


\section{Algebraic set theory}
\label{sec:ast}

The categorical set theories we have considered so far are analogues
of \zfc; in this section we consider categorical analogues of \nbg\
and \zfci.  One motivation for this is to find a more categorical way
to state the replacement axiom; another is to find an easier way to
deal with large categories.

We saw in \S\ref{sec:topos-indexed} that one appropriate notion of a
large category with respect to an elementary topos is an indexed
category.  As for na\"ive large categories in \zfc, the elementary
theory of indexed categories suffices when we only need to consider
them one at a time.  However, if we want to quantify over indexed
categories (such as in proving the Adjoint Functor Theorem), perform
constructions on them, or assemble them into a larger category (or a
2-category), we need to assume some sort of external set theory.  This
is automatic if our topos \bS\ is the category of small sets in \zfcs\
or \zfci---though in the latter case there is little need for indexed
categories.

A more categorical approach is to introduce a \emph{category of
  classes} \bC\ which contains the topos \bS\ of sets and also other
large categories.  This is a recent area of active research known as
\emph{algebraic set theory}; good introductions
are~\cite{awodey:outl-ast,jm:ast}.  The usual approach is to equip
\bC\ with a collection \cS\ of \emph{small morphisms}, the intuition
being that a morphism is small when it has small fibers.  Different
authors choose slightly different axioms on \bC\ and \cS, but in
general they can be classified as follows.
\begin{enumerate}
\item \bC\ has finite limits and colimits, and some additional level
  of structure allowing at least the interpretation of finitary logic.
\item Small maps are closed under composition, pullback, descent, and
  other basic constructions.
\item Every object $X$ of \bC\ has a \emph{powerclass} $\cP_s X$
  classifying its small subobjects.
\item There is a \emph{universal} object $U$ of \bC\ which `contains
  all small objects'.
\end{enumerate}
I will refer to some suitable version of these axioms collectively as
\textbf{Algebraic Set Theory (AST)}.  As with the elementary topos
axioms, the axioms of \ast\ can be augmented by well-pointedness,
choice, and existence of a \nno.

We define an object $X\in\bC$ to be \emph{small} if $X\to 1$ is a
small map.  Any model of \nbg\ gives rise to a model of \ast\ in which
the small objects are the sets.  Conversely, the axioms of \ast\ imply
that the category \bS\ of small objects is an elementary topos, which
inherits well-pointedness, choice, and a \nno\ from \bC.  Moreover, as
long as \bC\ is well-pointed, the \wptnc\ \bS\ also satisfies
replacement, in the sense described in \S\ref{sec:cst}.  We may not
quite get a model of \nbg\ by taking the objects of \bC\ to be the
classes, since some of them may be too large (violating limitation of
size), but in practice, working in \ast\ with well-pointedness,
choice, and a (small) \nno\ is essentially equivalent to working in
\nbg.

Note that just as in a topos, the internal logic of \bC\ is restricted
to bounded quantifiers.  However, now we can interpret `unbounded'
quantifiers ranging over all \emph{small} objects by using
\emph{bounded} quantifiers ranging over the universal object $U$.
This gives a categorical explanation of why comprehension in \nbg\ is
restricted to formulas that only quantify over sets.

If we assume in addition that \bC\ is itself a topos, as
in~\cite{streicher:universes}, then we obtain a theory equivalent to
what one might call \bzci.  Adding a replacement schema for \bC, as
described in \S\ref{sec:cst}, brings us up to \zfci.  In all cases, we
can define a \emph{large category} to be an internal category in \bC\
and a \emph{small category} to be one in \bS, and the development of
category theory then mirrors what happens in the corresponding
membership set theory.  (We could also use \bC\ to perform
constructions on \bS-indexed categories, but since \bS\ satisfies
replacement with respect to \bC\ there seems little need for the extra
complication.)


There are many other beautiful aspects to algebraic set theory, of
which I will only mention one: the cumulative hierarchy $V=\bigcup
V_\al$ and the class $\Omega$ of ordinals can be defined by universal
properties.  Define a \textbf{ZF-algebra} to be a partially ordered
class $A$ which has suprema of all \emph{small} families and is
equipped with a `successor' operation $s\maps A\to A$.  The class $V$
of all sets, ordered by inclusion and with $s(x) = \{x\}$, can then be
proven to be the initial ZF-algebra.  The class \Omega\ of ordinals is
also a ZF-algebra with $s(\al) = \al+1$, so we have a unique
homomorphism of ZF-algebras $\rho\maps V\to\Omega$; this turns out to
be essentially the \emph{rank} function.  The function $\Omega\to V$
sending \al\ to $V_\al$ can be characterized by an analogous universal
property, or as a right adjoint to \rho; see~\cite{jm:ast} for
details.


\begin{rmk}\label{rmk:constructive-2}
  So far I have focused exclusively on well-pointed categorical set
  theories with choice, because they are the most relevant to a
  mathematician looking for a categorical substitute for \zfc.  As
  noted in \autoref{rmk:constructive-1}, however, much of the
  independent interest of topos theory comes from the fact that
  \emph{any} elementary topos has an internal set-like logic, and in
  general this internal logic is not classical logic but
  \emph{constructive logic}.  Well-pointedness and choice are each
  quite special properties of a topos, and both independently imply
  that its logic is classical.

  Much of mathematics can be developed using constructive logic,
  although many classical definitions and results must be rephrased
  carefully to obtain a constructively meaningful or useful form.  For
  example, \emph{Tychonoff's theorem} that the product of compact
  spaces is compact is \emph{true} constructively, without the axiom
  of choice, but the definition of \emph{space} has to be modified;
  see~\cite[Part~C]{ptj:elephant}.  Classical concepts also often
  bifurcate into two or more inequivalent constructive ones.  For
  example, there are at least three different kinds of constructive
  ordinals with slightly different properties;
  see~\cite{jm:ast,taylor:ordinals}.

  Most relevantly for us, the axiom of replacement loses much of its
  power constructively: it no longer implies unbounded separation,
  Borel determinacy, or the usual sort of transfinite induction.
  Moreover, the categorical replacement axiom from
  \S\ref{sec:topos-indexed} is no longer sensible in the
  non-well-pointed case, and it is an open question whether it has
  some more general analogue.

  Algebraic set theory, like elementary topos theory, makes perfect
  sense (and is usually studied) without well-pointedness, but its
  version of replacement is also much weaker constructively.  In fact,
  \emph{any} elementary topos \bS\ can be embedded as the topos of
  small objects in a category \bC\ of classes, but the logic of \bC\
  will not in general be classical even if that of \bS\ is; thus every
  topos `satisfies replacement' in a constructive sense.  This is
  analogous to the use of indexed categories to `define away' the lack
  of replacement by considering only small families to begin with.
  The interested reader can learn more about constructive logic in
  toposes and \ast\ from the references cited above.
\end{rmk}

\section{Higher categories}
\label{sec:higher-cats}

One can envision more radically `categorical' foundations.  For
instance, it is hard to deny that in everyday mathematics we very
rarely care about large sets \emph{as sets}---we only care about them
insofar as they form the class of objects of some large
\emph{category}.  In particular, we only care about their elements up
to isomorphism.  Thus, in a sense, the `large' collections which arise
naturally in mathematics behave fundamentally differently than the
`small' ones.  This distinction is captured elegantly by the theory of
indexed categories, which gives an analogue of large categories
without making use of a prior notion of large set.

Analogously, we generally only care about large categories up to
equivalence, and thus we should regard them as objects of a 2-category
$\mathcal{CAT}$, rather than of a category $\mathbf{CAT}$.  This seems
to suggest that instead of axiomatizing the category of classes, a
more categorical generalization of elementary toposes would be to
axiomatize the 2-category of large categories.  Notable steps have
been made in this direction
(see~\cite{street:elem-cosmoi-i,weber:2toposes}), but I think it is
fair to say that a fully satisfactory answer has not yet emerged.
Since this approach also requires a good deal of familiarity with
2-categories, I will not attempt to explain it further here.

Of course, once we start to study 2-categories, we need to assemble
them into 3-categories, and so on \emph{ad infinitum}.  Philosophical
and mathematical remarks along these lines can be found
in~\cite{mak:catfdn}, among other places.


\section{Conclusion}
\label{sec:conclusion}

We have explored many possible foundations for category theory, including:
\begin{enumerate}[(1)]
\item A na\"ive approach which remains within \zfc.
\item Introducing classes as objects, as in \nbg\ and \mk.
\item Using an inaccessible to distinguish small sets from large ones (\zfci).
\item Using a reflection principle, perhaps combined with
  inaccessibles, as in \zfcs\ and \zmcs.
\item Categorical versions of the above, using toposes (\etcs) or
  categories of classes (\ast).
\end{enumerate}
Each has advantages and disadvantages, and I do not mean to put one
forward as \emph{the} correct foundation for category theory; I leave
that choice to the reader's aesthetic and mathematical judgment.

Instead, let me end by reiterating that for the \emph{basic} theorems
of category theory, the choice of foundation is essentially
irrelevant.  Each of the above proposals deals with the distinction
between small and large in a way which is fully satisfactory for
proving results such as the Adjoint Functor Theorem (except that in
some cases we have to state it as a meta-theorem, or add
small-definability restrictions).  However, as we have seen, the
choice of foundation does matter for some more elaborate
constructions.  Thus, I believe it is important for students and users
of category theory to have \emph{some} familiarity with its possible
foundations, and I hope to have partially addressed that need here.

\bibliographystyle{halpha}
\bibliography{all,shulman}

\end{document}